\documentclass[a4paper,11pt,reqno]{amsart} 
\usepackage{amssymb,amsthm,amsmath}
\usepackage{mathtools}
\usepackage{amsfonts}
\usepackage{amscd}
\usepackage{amssymb}

\usepackage{enumerate}
\usepackage{bbm}
\usepackage{graphicx}
\allowdisplaybreaks
\usepackage{mathabx}
\usepackage{color}
\usepackage{amsbsy}
\usepackage{graphicx}
\usepackage{amsthm}
\usepackage{amsmath}
\usepackage{amsxtra}
\usepackage{mathrsfs}
\usepackage{bbm}
\usepackage{dsfont}
\usepackage{enumitem}
\usepackage{comment}

\usepackage{enumerate}
\usepackage{xcolor}
\usepackage{float} % To use [H]

\usepackage{url}
\usepackage{soul}
\usepackage[hidelinks]{hyperref} 

\usepackage{ifthen}
\newtheorem{theorem}{Theorem}[section]
\newtheorem{defi}[theorem]{Definition}
\newtheorem{remark}[theorem]{Remark}
\newtheorem{corollary}[theorem]{Corollary}
\newtheorem{prop}[theorem]{Proposition}
\newtheorem{lemma}[theorem]{Lemma}

\numberwithin{equation}{section}
\definecolor{red}{rgb}{1.0, 0.0, 0.0}
\newcommand{\Bea}{\begin{eqnarray*}}
	\newcommand{\Eea}{\end{eqnarray*}}
\newcommand{\Be} {\begin{equation*}}
	\newcommand{\Ee} {\end{equation*}}
\newcommand{\be} {\begin{equation}}
	\newcommand{\ee} {\end{equation}}
\newcommand{\bea} {\begin{eqnarray}}
	\newcommand{\eea} {\end{eqnarray}}

\newcommand{\lab} {\label}

\newcommand{\abs}[1]{\left\vert#1\right\vert}

\newcommand{\al}{\alpha}

\newcommand{\la}{\lambda}

\title[ Some Versions of Beurling's Theorem on H-type Groups ]{Some Versions of Beurling's Theorem on H-type Groups}
\author[A. Dasgupta]{Aparajita Dasgupta}
\address{
	Aparajita Dasgupta:
	\endgraf
	Department of Mathematics
	\endgraf
	Indian Institute of Technology Delhi, Hauz Khas
	\endgraf
	New Delhi-110016 
	\endgraf
	India
	\endgraf
	{\it E-mail address:} {\rm adasgupta@maths.iitd.ac.in}
}
\author[P. Gulia]{Prerna Gulia}
\address{
	Prerna Gulia:
	\endgraf
	Department of Mathematics
	\endgraf
	Indian Institute of Technology Delhi, Hauz Khas
	\endgraf
	New Delhi-110016 
	\endgraf
	India
	\endgraf
	{\it E-mail address:} {\rm prernagulia64@gmail.com}
}
\author[S. Pusti]{Sanjoy Pusti}
\address{
	Sanjoy Pusti:
	\endgraf
	Department of Mathematics
	\endgraf
	Indian Institute of Technology Bombay
	\endgraf
	Mumbai-400076 
	\endgraf
	India
	\endgraf
	{\it E-mail address:} {\rm sanjoy@math.iitb.ac.in}
}
\author[S. Thangavelu]{Sundaram Thangavelu}
\address{
	Sundaram Thangavelu:
	\endgraf
	Department of Mathematics
	\endgraf
	Indian Institute of Science, Bangalore
	\endgraf
	Bangalore-560012 
	\endgraf
	India
	\endgraf
	{\it E-mail address:} {\rm veluma@iisc.ac.in}
}

\vspace{1cm}

\allowdisplaybreaks
\begin{document}
	
	\begin{abstract}
We prove an analogue of Beurling’s theorem on the H-type groups of certain dimensions after establishing the Gutzmer’s formula for the H-type groups. We also obtain some other versions of the theorem using the modified Radon transform.

 % In 2007, Thangavelu proved on the Heisenberg group, the Gutzmer’s formula, which was first established by Faraut for all non-compact Riemannian symmetric spaces. This formula was used to establish an analogue of the Beurling’s theorem on the Heisenberg group. In this article, we extend the Gutzmer’s formula to H-type groups, notion introduced by Kaplan. We obtain the formula by defining H-type motion group and a family of its class 1 representations. As an application,	% 
 
 \end{abstract}
	\maketitle
%	\tableofcontents

		\section{Introduction}
The uncertainty principle states that a function and its Fourier transform cannot both decay very rapidly at the same time. The most general version of the uncertainty principles is Beurling's theorem. Beurling’s theorem on \(\mathbb{R}^n\) states that there is no nontrivial function which satisfies
\Be
 \int_{\mathbb{R}^n}{\int_{\mathbb{R}^n} \lvert f(y) \rvert \lvert \Hat{f}(\xi) \rvert e^{\lvert (y,\xi) \rvert}}\;dy\; d\xi < \infty.
\Ee
Several other uncertainty principles  such
 as the theorems of Hardy, Morgan, Cowling and Price, Gelfand and Shilov follows from this Beurling's theorem (\cite{sarkarbeurling}).
%To extend Beurling's Theorem on the Heisenberg group \cite{thangavelu-beurling}, the author imposed the following condition on the Fourier transform of the function \(f\) which helps in holomorphic extension of the function \(f^{\la}\), where, \(f^{\la}\) is the inverse Fourier transform of \(f\) in the center variable 
% \begin{equation}\label{eqn1.1}
% \int_K {\lVert \pi_\la(\sigma(z,w))\Hat{f}(\la) \rVert}^2_{HS} \;d\sigma\;<\infty 
 %\end{equation}
 
 For a linear operator \(T\) on \(L^2(\mathbb{R}^n)\), let
 \(\mathcal{D}(T)\) and \(\mathcal{R}(T)\) stand
 for its domain and range. Let\({\lVert T \rVert}_1 =tr\lvert T\rvert \) stand for the
 trace norm of the operator \(T\) where \(\lvert T\rvert =(T^*T)^{1/2}\). Also, let \(d\sigma\) be the Haar measure on \(K=Sp(n,R) \cap O(2n,R) \cong U(n)\). The Beurling theorem on \(\mathbb{H}^n\) (proved by Thangavelu in \cite{thangavelubeurling}) states that:
 \begin{theorem}
     [\textbf{Beurling's theorem on \(\mathbb{H}^n\)
    }]\label{beurling on heisenberg}
     For $f \in L^2\left(\mathbb{H}^n\right)$, assume that $\mathcal{R}(\hat{f}(\lambda)) \subset \mathcal{D}\left(\pi_\lambda(z, w)\right)$ and $\pi_\lambda(z, w) \hat{f}(\lambda)$ is of trace class for all $(z, w) \in \mathbb{C}^{2 n}$ such that
   \begin{equation}\label{eqn1.1}
 \int_K {\lVert \pi_\la(\sigma(z,w))\Hat{f}(\la) \rVert}^2_{HS} \;d\sigma\;<\infty.
 \end{equation}  
      If the condition
$$
\int_{\mathbb{R}^{2 n}}\left|f^\lambda(y, v)\right| \quad\left(\sup _{\left|\left(y^{\prime}, v^{\prime}\right)\right| \leq|(y, v)|}\left\|\pi_\lambda\left(i\left(y^{\prime}, v^{\prime}\right)\right) \hat{f}(\lambda)\right\|_1\right) d y d v<\infty,
$$
holds for almost every $\lambda \in \mathbb{R}^*$, then $f=0$.
\end{theorem}
\noindent Later, this was used to establish Hardy, Cowling and Price theorems on the Heisenberg group. Moreover, the above theorem is one of the several versions of the Beurling's theorem that have been established on the Heisenberg group. In this paper, we will also reflect upon one more version of the same, which was established by Parui and Sarkar (\cite{Beurling-rudra-parui}) in 2008. The authors proved the following modified version of the theorem on the Heisenberg group.
\begin{theorem}[\textbf{Modified Beurling's theorem on \(\mathbb{H}^n\)
    }]\label{Beurling-Rudra-Parui}
  Let \(f\in L^2(\mathbb{H}^n)\) and for some \(M_1,M_2\geq 0\), the function \(f\) satisfies
  $$
 \int\limits_{\mathbb{H}^n}{\int\limits_{\mathbb{R}}\frac{|f(z,a)|{\lVert \Hat{f}(\la)\rVert}_{HS}e^{\lvert a \rvert \lvert  \la|}{\lvert \la \rvert}^n}{{(1+|z|)}^{M_1}{(1+|\la|+|a|)}^{M_2}}d\la dz d a}< \infty. 
  $$
  Then \(f(z,a)=e^{-c{a}^2}{(1+|z|)}^{M_1}\left (\sum\limits_{j=0}^{l}\psi_j(z){a}^j\right )\), where \(c\) is a constant, the functions \(\psi_j \in L^1(\mathbb{C}^n)\cap L^2(\mathbb{C}^n)\) and \(l<\frac{M_2-n/2-1}{2}.\)
\end{theorem} 

Our aim in this paper is to prove an analogue of Theorem \ref{beurling on heisenberg} and Theorem \ref{Beurling-Rudra-Parui} on an H-type group $N$. The class of H-type groups has received significant attention in the area of harmonic analysis. Kaplan and Ricci \cite{kaplanricci} established the framework for the study of the representation theory of H-type groups, drawing inspiration from the Bargmann representations of the Heisenberg group. 
%They also raised the question of when does the pair \((A(N),N)\), where \(A(N)\) is the group of automorphisms of H-type group \(N\) that act as orthogonal transformations on its Lie algebra, forms a Gelfand pair. Both positive and negative results were obtained in this context depending on the degree and irreducibility of the Lie algebra. In 1985, Ricci \cite{riccicomm} further pursued this to classify the groups of H-type for which the algebra of \(L^1\) functions, that are \(A(N)\) invariant, is commutative. 
%In the last two decades, the heat kernel \cite{yangheatkernel}, heat kernel on nilmanifolds associated to H-type groups \cite{dasgupta-thangavelu-nilmanifold}, restriction theorem \cite{LiurestrictionHtype}, spherical means \cite{narayananinjectivity,yasser2024gelfand} and Weyl transform \cite{dasguptawong} on the H-type group has been extensively studied.

For an H-type group $N$, let $\mathfrak{n}$ be it's Lie algebra and $\mathfrak{z}$ be the center of $\mathfrak{n}$ with $m=\text{dim}\,\mathfrak{z}$ and $\text{dim }\mathfrak{z}^{\perp}=2n$. For an unit vector \(\omega \in \mathfrak{z}\), let \(k(\omega)\) denote the orthogonal
 complement of \(\omega\) in \(\mathfrak{z}\). Let \(q_{\omega}\) denote the quotient map from \(\mathfrak{n}\) to \(\mathfrak{n}/k(\omega)\). It is known that this quotient algebra can be identified with \(\mathfrak{v}\oplus \mathbb{R}\) and is isomorphic to the Heisenberg Lie algebra \(\mathfrak{h}_{\omega}\). We denote the associated Lie group by \(\mathbb{H}^n_{\omega}\). Let \(A(N)\) denote the group of automorphisms of \(N\) that act as orthogonal transformations on the Lie algebra \(\mathfrak{n}\). 
Let \(U\) denote the subgroup of automorphisms that act trivially on the center, that is, 
$$
U=\{\theta \in A(N): \theta(z)=z,\; \forall z \in \mathfrak{z}\}.
$$ For the case of Heisenberg group $\mathbb H^n$, this $U$ turns out to be $\mathrm{U}(n)$. For each $\boldsymbol{\lambda}\neq0$ in $\mathbb R^m$, let $\pi_{\boldsymbol{\lambda}}$ be the Schr\"{o}dinger representation of the group $N$ on $L^2(\mathbb R^n)$. 
This $\pi_{\boldsymbol{\lambda}}$ can be written as \(\pi_{\la,\omega}\) for some \(\la\in \mathbb R^\ast, \omega \in \mathbb{S}^{m-1}\) (see \cite{roncal-thangavelu-holomorphic-extension} for details). Also, the restriction of \(\pi_{\la,\omega}\) to \(\mathbb{H}^n_{\omega}\)  is unitarily equivalent to the Schr\"odinger representation of the Heisenberg group.

In this paper, we first prove the following analogue of Theorem \ref{beurling on heisenberg} for the H type group $N$.
\begin{theorem}\label{first Beurling's using radon transform}
Let \(f\in L^2(N)\) such that \(\mathcal{R}(\Hat{f}(\la,\omega)) \subset \mathcal{D}(\pi_{\la,\omega}(z,\gamma))\) and \(\pi_{\la,\omega}(z, \gamma)\Hat{f}(\la,\omega)\) is of trace class for all $(z, \gamma)\in \mathbb C^{2n}$, such that for each \(\omega \in \mathbb{S}^{m-1}\),
\begin{equation}\label{cond1}
   \displaystyle\int_{U(n)}{\left \lVert \pi_{\la,\omega}\left({q_{\omega}}^{-1}(\sigma(z,\gamma)\right)\Hat{f}(\la,\omega)\right \rVert}^2_{HS} \,d\sigma< \infty.
\end{equation} 
If for each \(\omega \in \mathbb{S}^{m-1}\), the function satisfies 
    \begin{equation}\label{cond2}
      \int_{\mathbb{R}^{2n}}\left | f^{\la\omega}(y,v)  \right | \left ( \sup_{|(y',v')| \leq |(y,v)|} {\left \lVert \pi_{\la,\omega}(i(y',v'))\Hat{f}(\la,\omega)\right \rVert }_1  \right ) dydv < \infty,  
    \end{equation}
    for almost every \(\la \in \mathbb{R}^*\), then the function \(f=0\).
\end{theorem}

The condition (\ref{cond1}) in the theorem above is not natural as we integrated on $\mathrm{U}(n)$. The integration should have been on the group $U$, not on $\mathrm{U}(n)$. We had to take this as we use Radon transform and transferred the problem to the Heisenberg group. 

We got the exact version of the Beurling's theorem in the following theorem for the case when $m=2, 3$. 

\begin{theorem}\label{Beurling's theorem}
    Let \(f\) be a Schwartz class function on \(N\) such that  \(\mathcal{R}(\Hat{f}(\boldsymbol{\lambda})) \subseteq \mathcal{D}(\pi_{\boldsymbol{\lambda}}(z,w))\) and \(\pi_{\boldsymbol{\lambda}}(z,w)\Hat{f}(\boldsymbol{\lambda}) \text{ is a trace class operator for all } (z,w) \in \mathbb{C}^{2n}\) such that, \Be \displaystyle\int_{U}{\lVert \pi_{\boldsymbol{\lambda}}(\sigma(z,w))\Hat{f}(\boldsymbol{\lambda})\rVert}^2_{HS} \,d\sigma< \infty. \Ee If 
    $$
    \int_{\mathbb{R}^{2n}}\left | f^{\boldsymbol{\lambda}}(y,v)  \right | \left ( \sup_{|(y',v')| \leq |(y,v)|} {\left \lVert \pi_{\boldsymbol{\lambda}}(i(y',v'))\Hat{f}(\boldsymbol{\lambda})\right \rVert }_1  \right ) dydv < \infty,
    $$
    for almost every \(\boldsymbol{\lambda} \in \mathbb{R}^m\setminus \{0\}\) then the function \(f=0\).
\end{theorem}
To prove the theorem above, we use Gutzmer's formula and ideas of Thangavelu \cite{thangavelubeurling}.  Roughly speaking, Gutzmer's formula is the Plancherel theorem for holomorphic functions. The Gutzmer formula for the Heisenberg group has already been established (\cite{gutzmer}). 

In general, when dealing with a H-type group \( N \), the analysis becomes significantly more involved due to the complexity of their automorphism groups.  The automorphism group involves the spin group, which is more complex in nature.
A crucial component in establishing the Gutzmer's formula is a specific estimate involving rank one operators (see \cite[Proposition 4.3]{gutzmer}).  One of the defining features of a Gelfand pair \( (K,G) \) is that for every irreducible representation \( \pi \) of \( G \), the subspace of \( K \)-invariant vectors has dimension at most one. In particular, if \( \pi \) is a class one representation, it contains (up to scalar multiple) a unique \( K \)-invariant vector, often referred to as the fixed vector. This property is instrumental in proving that the Fourier transform of a Schwartz class function is a rank one operator.

In the case of the Heisenberg group, Thangavelu considers the subgroup \(U=U(n)\) of the automorphism group (instead of whole automorphism group of Heisenberg group). This is relatively simple to handle. % in \cite{gutzmer} expressed The Fourier inversion formula for the Heisenberg group was expressed in terms of certain class one representations of the Heisenberg motion group \(U(n) \ltimes N\), the semi-direct product of the Heisenberg group with the unitary group \( U(n) \). 
Similarly, to avoid complexity (involving spin group) we consider the pair \((U,N)\) (instead of $( A(N),N)$) which is Gelfand if and only if \(m=1,2,3\). Thus, we restrict ourselves to these three cases only (for Gutzmer's formula and subsequently for Theorem \ref{Beurling's theorem}). We first establish Gutzmer’s formula in the case when \( m = 3 \), and subsequently for the case \( m = 2 \). These two cases are treated separately due to the different behavior of the action of the compact group \( U \) involved in the analysis. Specifically:

\begin{itemize}
    \item For \( m = 3 \), the group \( U \) acts transitively on the product of spheres, which plays a critical role in the decomposition of functions and in the application of harmonic analysis. 
    
    \item In contrast, for \( m = 2 \), the group \( U \) acts transitively on spheres centered at the origin.
\end{itemize}
The case $m=1$ gives the pair $(\mathbb H^n, U(n))$ which is dealt in (\cite{gutzmer, thangavelubeurling}).

Our next theorem is an analogue of Theorem \ref{Beurling-Rudra-Parui}. We prove this theorem again by using Radon transform, but unfortunately, we could not recover the function, instead we got the Radon transform of the function. 
\begin{theorem}\label{second Beurling's using radon transform}
Let \(f\in L^1\cap L^2(N)\) such that for some \(M_1>n ,M_2 \geq 0\) and each \(\omega \in \mathbb{S}^{m-1}\), 
\begin{equation}\label{second Beurling's eq1}
\int\limits_{\mathbb{R}^{2n}\times \mathbb{R}^m}{\int\limits_{\mathbb{R}}\frac{\left \lvert (-\Delta_t)^{(m-1)/4}f(z,t)\right \rvert{\lVert \Hat{f}(\la,\omega)\rVert}_{HS}e^{\lvert t \rvert \lvert  \la|}{\lvert \la \rvert}^{n+(m-1)/2}}{{(1+|z|)}^{M_1}{(1+|\la|+|t|)}^{M_2}}d\la dz dt}< \infty. 
\end{equation}
Then for each \(\omega \in \mathbb{S}^{m-1}\), \Be R_{\omega}f(z,a)=e^{-c_{\omega}{a}^2}{(1+|z|)}^{M_1}\left (\sum\limits_{j=0}^{l}\psi_{j,\omega}(z){a}^j\right ), \Ee where \(c_{\omega} \) is a constant, the functions \(\psi_{j,\omega} \in L^1(\mathbb{C}^n)\cap L^2(\mathbb{C}^n)\) and \(l<\frac{M_2-n/2-1}{2}.\)
\end{theorem}

We prove Theorem \ref{first Beurling's using radon transform} and Theorem \ref{second Beurling's using radon transform} using Radon transform.  The Radon transform was first introduced by Strichartz \cite{strichartz} in the study of the Heisenberg group. The strategy is to use the Radon transform in the central variable of the H-type group $N$ to reduce the problem to the Heisenberg group.
 
 We now conclude this introduction by describing the paper's outline. In Section \ref{sec2}, we present essential background information on H-type groups, encompassing their group operation, representation theory, and Fourier transform. In section \ref{sec3}, we establish two versions of Beurling's theorem using Radon transform, namely Theorem \ref{first Beurling's using radon transform} and Theorem \ref{second Beurling's using radon transform}. Section \ref{sec4} focuses on constructing the H-type motion group and a family of its class one representations, leading to the establishment of Gutzmer's formula. In Section \ref{sec5}, we introduce the reduced H-type group and provide the necessary estimates to support the Beurling theorem. Finally, in Section \ref{sec6}, we prove Beurling's theorem using Gutzmer's formula.
  
  \section{Preliminaries}\label{sec2}
 In this section, we present the essential preliminaries for H-type groups including Fourier transform and the automorphism groups of H-type groups. These have been covered in detail in \cite{bonfi,dasguptawong,functionalhtype,yangheatkernel, mingkai}.
  \begin{defi}
      A Heisenberg type or H-type algebra is a finite-dimensional real Lie algebra \(\mathfrak{n}\), equipped with an inner product $\left \langle \cdot,\cdot\right \rangle$ such that $[\mathfrak{z}^\perp, \mathfrak{z}^\perp]=\mathfrak{z}$, where $\mathfrak{z}$ denotes the center of $\mathfrak{n}$ and for every $z \in \mathfrak{z}$, the map $J_z: \mathfrak{z}^\perp \to \mathfrak{z}^\perp$ defined by 
   \begin{center}
       $\left \langle J_z (v),w \right \rangle=\left \langle z, [v,w] \right \rangle$
   \end{center}
   is orthogonal whenever $ \lvert z \rvert=1$.
  \end{defi}
 %Throughout this paper, \(|\la|\) refers to the norm of a vector \(x\) induced by the inner product of the corresponding inner product space.
Throughout this paper we will be using notation \(v \cdot w\) for inner product of two vectors \(v\) and \(w\).\\
A connected and simply connected Lie group $N$ is said to be of H-type if its Lie algebra is of H-type. Let \(N\) be an H-type group and \(\mathfrak{n}\) be its corresponding Lie algebra. Then \(\mathfrak{n}\) can be written as orthogonal sum of \(\mathfrak{z}\text{ and }\mathfrak{z}^\perp\), i.e., \(\mathfrak{n}=\mathfrak{z}+\mathfrak{z}^{\perp}\), the center and its orthogonal complement. Furthermore, it is known that dim \(\mathfrak{z}^\perp\) has to be even, say \(2n\). Thus, we identify \(\mathfrak{z}^\perp\) with \(\mathbb{C}^n\) or \(\mathbb{R}^{2n}\) and \(\mathfrak{z}\) with \(\mathbb{R}^m\).% We write \((z,t)\) for points in \(N\) and the group law of \(N\) can be written as 

 We identify the Lie group \(N\) with its nilpotent
 Lie algebra \(\mathfrak{n}\) via the exponential map and write \((z,t)\) for points in \(N\), where \(z \in \mathbb{C}^n\) and
 \(t \in \mathbb{R}^m\). It follows from the Baker-Campbell-Hausdorff formula that the group law on \(N\) is given by
 \begin{equation}\label{group law}
    (z,t)(w,s)=\left(z+w,t+s+\frac{1}{2} [z,w] \right),
\end{equation}
where \([z,w]_i=\left \langle z, U^i w \right \rangle\) (see \cite{bonfi}), for matrices \(U^i\) satisfying the following conditions:
\begin{enumerate}
    \item \(U^i\) is a skew-symmetric orthogonal matrix for \(1 \leq i \leq m\),
    \item \(U^iU^j+U^jU^i=0\) for \(i \neq j\).
\end{enumerate}
The Haar measure on \(N\) is the Lebesgue measure on \(\mathfrak{n}\) and is denoted
 by \(dzdt\).
It is important to note that we can identify \(\mathfrak{n}\) with \(\mathbb{C}^n \times \mathbb{R}^m\) via any orthonormal basis of \(\mathfrak{n}\). However, for some computations, we require a certain choice of orthonormal basis (see \cite{functionalhtype,narayananinjectivity}) to be described later.
The convolution between two functions $ f, g \in L^1(N) $ is defined as
$$ f \ast g(z,t) = \int_{\mathbb{C}^n} \int_{\mathbb{R}^m} f(z-w, t-s-\frac{1}{2}[z,w])\, g(w,s) \, dw\, ds.$$
For any $ \boldsymbol{\lambda} \in \mathbb{R}^m \setminus \{0\}$, let $ f^{\boldsymbol{\lambda}}(z) $ stand for the inverse Fourier transform of $ f $ in the central variable. Thus 
$$ f^{\boldsymbol{\lambda}}(z) = \int_{\mathbb{R}^m} f(z,t) \, e^{i \boldsymbol{\lambda} \cdot t} \, dt.$$
An easy calculation shows that 
$$ (f \ast g)^{\boldsymbol{\lambda}} (z) = \int_{\mathbb{C}^n} f^{\boldsymbol{\lambda}}(z-w)\, g^{\boldsymbol{\lambda}}(w)\, e^{\frac{i}{2}  \boldsymbol{\lambda}\cdot [z,w]}\, dw.$$ 
We will refer to the integral on the right-hand side as the $ \boldsymbol{\lambda}$-twisted convolution of $ f^{\boldsymbol{\lambda}} $ with $ g^{\boldsymbol{\lambda}}$, denoted by $ f^{\boldsymbol{\lambda}} \times_{\boldsymbol{\lambda}} g^{\boldsymbol{\lambda}}(z).$ Note that here \(\boldsymbol{\lambda} \in \mathbb{R}^m\), unlike the usual twisted covolution where \(\lambda \in \mathbb{R}\). The Lie bracket in the above expression is determined by the H-type group \(N\), but we will see below that with a suitable choice of basis,  \(\boldsymbol{\lambda}\)-twisted convolution takes a simple form as in Eq. \eqref{transformed twisted convolution}.

It is known that for each $ \boldsymbol{\lambda} \in \mathbb{R}^m, \boldsymbol{\lambda} \neq 0,$ there is an irreducible unitary representation $ \pi_{\boldsymbol{\lambda}} $ of $ N $ realised on $ L^2(\mathbb{R}^n).$  The group Fourier transform of a function $ f \in L^1(N) $ is then defined to be the operator valued function
$$\pi_{\boldsymbol{\lambda}}(f)= \widehat{f}(\boldsymbol{\lambda}) = \int_{\mathbb{C}^n}\int_{\mathbb{R}^m} f(z,t)\, \pi_{\boldsymbol{\lambda}}(z,t)\, dz\, dt.$$
Since $ \pi_{\boldsymbol{\lambda}}(z,t) = \pi_{\boldsymbol{\lambda}}(z,0)\pi_{\boldsymbol{\lambda}}(0,t) = e^{i \boldsymbol{\lambda} \cdot t}\pi_{\boldsymbol{\lambda}}(z,0)$, it follows that 
$$  \widehat{f}(\boldsymbol{\lambda}) = \int_{\mathbb{C}^n} f^{\boldsymbol{\lambda}}(z)\, \pi_{\boldsymbol{\lambda}}(z,0)\, dz =: W_{\boldsymbol{\lambda}}(f^{\boldsymbol{\lambda}}).$$ 
The operator $ W_{\boldsymbol{\lambda}}(f^{\boldsymbol{\lambda}}) $ is known as the Weyl transform of $ f^{\boldsymbol{\lambda}}.$  It is then easy to check that
\begin{equation*}
W_{\boldsymbol{\lambda}} (\varphi \times_{\boldsymbol{\lambda}} \psi) = W_{\boldsymbol{\lambda}}(\varphi) W_{\boldsymbol{\lambda}}(\psi),
\end{equation*}  for any $ \varphi, \psi \in L^1(\mathbb{C}^n).$ In order to describe the representation $ \pi_{\boldsymbol{\lambda}} $ and the operator $ W_{\boldsymbol{\lambda}} $ more explicitly we now choose a special orthonormal basis for  \(\mathfrak{z}\text{ and }\mathfrak{z}^\perp\).
For each \(0 \neq \boldsymbol{\lambda} \in \mathbb{R}^m\), consider the skew symmetric bilinear form on \(\mathfrak{z}^\perp\) defined by the relation
$$
\left \langle J_{\boldsymbol{\lambda}} u,v \right \rangle =\boldsymbol{\lambda}([u,v]).
$$
Additionally, choose an orthonormal basis \(\{E_1(\boldsymbol{\lambda}), \dots, E_n(\boldsymbol{\lambda}),E'_1(\boldsymbol{\lambda}), \dots, E'_n(\boldsymbol{\lambda})\}\) of \(\mathfrak{z}^\perp\) such that 
\begin{align}
J_{\boldsymbol{\lambda}} E_j(\boldsymbol{\lambda})&=-E'_j(\boldsymbol{\lambda})\nonumber,\\
J_{\boldsymbol{\lambda}} E'_j(\boldsymbol{\lambda})&= E_j(\boldsymbol{\lambda}). \label{basis}
\end{align}
Also, choose an orthonormal basis \(\{ \boldsymbol{\lambda}= e_1(\boldsymbol{\lambda}),\dots, e_m(\boldsymbol{\lambda})  \}\) of \(\mathfrak{z}\), such that \(\boldsymbol{\lambda} \cdot e_1(\boldsymbol{\lambda})=|\boldsymbol{\lambda}|\) and \(\boldsymbol{\lambda} \cdot e_j(\boldsymbol{\lambda})=0\) for \(j \neq 1\). 
If \(\mathfrak{n}\) is identified with \(\mathbb{C}^n \times \mathbb{R}^m\) with the orthonormal basis chosen as above, the Lie bracket satisfies
\begin{equation*}
[z,w]_1=\langle z,U^1 w \rangle=\operatorname{Im}(z \cdot \overline{w}).\,\, %\textcolor{red}{[z,w]_j =0, j\neq 1}.
\end{equation*}
%\textcolor{red}{If we identify $ N $ with $ \mathbb{C}^n \times \mathbb{R}^m$ using the above coordinate system, then  the group law takes the form 
%$$ (z,a) (w,b) = (z+w, a+b+ \frac{1}{2} \operatorname{Im}(z \cdot \overline{w})\omega).$$}
%For any fixed $ \omega $ if we let $ N_\omega = \mathbb{C}^n \times \mathbb{R} \omega$ then it is clear that $ N_\omega $ is isomorphic to the standard Heisenberg group $ \mathbb H^n_{\omega}.$
In the above coordinate system, the definition of the $ \boldsymbol{\lambda}$-twisted convolution takes the form 
\begin{equation}\label{transformed twisted convolution}
 \varphi \times_{\boldsymbol{\lambda}} \psi(z) = \int_{\mathbb{C}^n} \varphi(z-w)\, \psi(w)\, e^{\frac{i}{2} |\boldsymbol{\lambda}| \operatorname{Im}(z \cdot \overline{w})}\, dw.   
\end{equation} 
When $ \omega $ is fixed and $ \boldsymbol{\lambda} = \la \omega,$ we suppress $ \omega $ and simply write $  \varphi \times_\la \psi$ instead of $  \varphi \times_{\la\omega} \psi$. 
With this understanding, we have the relation 
$$ (f \ast g)^{\la \omega}(z) = f^{\la \omega} \times_\la g^{\la\omega}(z).$$
We will now describe a family of irreducible unitary representations of $ N$. For any non-zero $ \la \in \mathbb{R} $, let $ \pi_\la$ be the Schr\"odinger representation of the Heisenberg group $ \mathbb H^n $ defined by 
$$ \pi_\la(z,a) \varphi(\xi) = e^{i\la a}\, e^{i\la ( x \cdot \xi+\frac{1}{2} x \cdot y)} \varphi(\xi+y),\, \varphi \in L^2(\mathbb{R}^n).$$ 
This $ \pi_\la $ gives rise to the representation $ \pi_{\la, \omega} $ of $ N $. Using the coordinate system described above, we define
$$ \pi_{\la, \omega}(z,t) \varphi(\xi) = e^{i\la t \cdot \omega}\, e^{i\la ( x \cdot \xi+\frac{1}{2} x \cdot y)} \varphi(\xi+y),\, \varphi \in L^2(\mathbb{R}^n).$$
It is then easy to check that $ \pi_{\la, \omega} $ is indeed an irreducible unitary representation of $ N $ realised on $ L^2(\mathbb{R}^n).$

%For each \(0 \neq \lambda \in \mathfrak{z}^*\) (identified with \(\mathbb{R}^m\)), consider the skew symmetric bilinear form on \(\mathfrak{z}^\perp\) defined by the relation
%$$
%\left \langle B(\lambda)u,v \right \rangle =\lambda([u,w]).
%$$
%Additionally, choose an orthonormal basis \(\{E_1(\lambda), \dots, E_n(\lambda),E'_1(\lambda), \dots, E'_n(\lambda)\}\) of \(\mathfrak{z}^\perp\) such that 
%\begin{align}
%B(\lambda)(E_i(\lambda))&=-|\lambda|E'_i(\lambda)\nonumber,\\
%B(\lambda)(E'_i(\lambda))&=|\lambda|E_i(\lambda). \label{basis}
%\end{align}
%Also, choose an orthonormal basis \(\{e_1,\dots, e_m  \}\) of \(\mathfrak{z}\) such that \(\lambda(e_1)=|\lambda|\) and \(\lambda(e_j)=0\) for \(j \neq 1\).
%If \(\mathfrak{n}\) is identified with \(\mathbb{C}^n \times \mathbb{R}^m\) with the orthonormal basis chosen as above, the Lie bracket satisfies
%\begin{equation*}
%[z,w]_1=\langle z,U^1 w \rangle=\operatorname{Im}(z \cdot \overline{w}).
%\end{equation*}
Let \(z=x+iy,\; x,y \in \mathbb{R}^n\). The left-invariant vector field of the algebra \(\mathfrak{n}\) of \(N\) that agree with \(\frac{\partial}{\partial x_i}\) and \(\frac{\partial}{\partial y_i}\) respectively at the origin are given by 
\begin{align*}
    X_i&=\frac{\partial}{\partial x_i}+\frac{1}{2}\sum\limits_{k=1}^{m}{\sum\limits_{l=1}^{2n}z_lU^k_{l,i}\frac{\partial}{\partial t_k}},\\
    Y_i&=\frac{\partial}{\partial y_i}+\frac{1}{2}\sum\limits_{k=1}^{m}{\sum\limits_{l=1}^{2n}z_lU^k_{l,i+n}\frac{\partial}{\partial t_k}}
\end{align*}
where \(z_l=x_l,\;z_{l+n}=y_l\) for \(1 \leq l \leq n\).\\
Then sublaplacian on the group \(N\) is given by 
$$
\mathcal{L}=-\sum\limits_{i=1}^{n}(X_i^2+Y_i^2)=-\Delta_z+\frac{1}{4}|z|^2T-\sum\limits_{k=1}^{m}\left \langle z,U^k \nabla_z \right \rangle T_k,
$$
where \Be
\Delta_z=\sum\limits_{i=1}^{2n}\frac{\partial^2}{\partial  z_i^2}, T_k=\frac{\partial}{\partial t_k}, T=-\left ( \sum\limits_{1}^{m}\frac{\partial^2}{\partial t_i^2}\right) \text{ and } \nabla_z={\left ( \frac{\partial}{\partial z_1}, \dots , \frac{\partial}{\partial z_{2n}}\right)}^t.
\Ee
%Once again, the above formulas will take a simpler form if we fix a unit vector $ \omega $ and choose an orthonormal basis as before.

In particular, one obtains the following lemma concerning the action of the sublaplacian \(\mathcal{L}\) on functions of the form \(e^{i \la \omega \cdot t} \phi(z)\). For a proof, we refer the reader to \cite[Lemma 1]{LiurestrictionHtype}.
%We denote \(|\la|= \sqrt{{\la_1}^2+ \cdots +{\la_m}^2}\).
\begin{lemma}\label{twisted laplacian lemma}
Given a vector $ \boldsymbol{\lambda}=\lambda\omega \in \mathfrak z$, with \(|\boldsymbol{\lambda}|=\lambda\) choose a coordinate system as before. If \(f(z,t) = e^{i \la \omega \cdot t} \phi(z)\), then
$ \mathcal{L}f(z,t) = e^{i \la \omega \cdot t}\mathcal{L}_{\lambda}\phi(z)$
 where,
$$
\mathcal{L}_{\la}=-\Delta_z+\frac{1}{4}\la^2{|z|}^2-i\la \sum_{j=1}^{n}\left( x_j\frac{\partial}{\partial y_j}-y_j\frac{\partial}{\partial x_j}\right)
$$
is the twisted Laplacian (or the scaled Hermite operator) on \(\mathbb{C}^{n}\).
\end{lemma}
%\textcolor{orange}{The Schr\"odinger representations of \(N\) \cite{yangheatkernel,mingkai} are parameterized by non zero elements \(\lambda\) in \(\mathbb{R}^m\) and are defined by 
%$$
%\pi_\la(x,u,t)\varphi(\xi)=e^{i\la \cdot t}e^{i |\la|(x \cdot\xi + \frac{1}{2}x\cdot u)} \varphi(\xi + u),\;\; \text{where } (x,u,t)\in N \text{ and } \varphi \in L^2(\mathbb{R}^n).
%$$
%Observe that \(\pi_\la(0,0,t)=e^{i \la \cdot t}I\). Moreover, by Stone Von-Neumann theorem, any irreducible, unitary representation which is equal to \(e^{i \la \cdot t} I\), on the center of \(N\), is unitarily equivalent to \(\pi_\la\) (\cite{dasguptawong}).
%\begin{remark}
 %   We note that \(\pi_\la(x,u,0)=\rho_{|\la|}(x,u,0)\) where \(\rho_{|\la|}\) refers to the Schr\"odinger representation of the Heisenberg group.
%\end{remark}
%Let \(\la \in {\mathbb{R}^m}^*\) and \(f \in L^1(N)\). The Fourier transform \(\Hat{f}(\la)\) of the function \(f\), is the bounded operator on \(L^2(\mathbb{R}^n)\) defined as 
%$$
%\pi_\lambda(f)=\Hat{f}(\lambda) = \int\limits_{\mathbb{R}^{2n+m}} f(z,t)\;\pi_{\lambda}(z,t) dz dt.
%$$}
We have the following analogue of Plancherel formula (\cite{dasguptawong}) on \(N\): For $f \in L^2(N)$,  \(\Hat{f}(\boldsymbol{\lambda})\) is a Hilbert–Schmidt
 operator on $L^2(\mathbb{R}^n)$ for all \(0 \neq \boldsymbol{\lambda} \in \mathbb{R}^m\) with
 \begin{center}
    $ {\lVert f \rVert}^2_2 = \displaystyle\int\limits_{\mathbb{R}^m} {\lVert \Hat{f}(\boldsymbol{\lambda})\rVert}^2_{HS}\; d\mu(\boldsymbol{\lambda})$,
 \end{center}
 where $d\mu(\boldsymbol{\lambda})={\lvert \boldsymbol{\lambda} \rvert}^n (2\pi)^{-(n+m)} d\boldsymbol{\lambda} $ is the Plancherel measure.
Note that in terms of representations \(\pi_{\lambda,\omega}\), \(\lambda>0, \omega \in \mathbb{S}^{m-1}\) we have the following Plancherel formula
$$
\| f\|_2^2 = (2\pi)^{-(n+m)}\int_{0}^\infty \int_{S^{m-1}}  {\lVert \Hat{f}(\lambda, \omega)\rVert}^2_{HS}\;  \, |\lambda|^{n+m-1}\,d\la\, d\sigma(\omega).
$$
For \(\al,\beta \in \mathbb{N}^n\), let \(\Phi_\al\), \(\Phi_{\al, \beta}\) be the normalised Hermite function and the special Hermite function on \(\mathbb{R}^n\) respectively (\cite[Page 18]{thangaveluheisenberg}).
For $\boldsymbol{\lambda}=\la \omega,\; |\boldsymbol{\lambda}|=\lambda$, we define 
$$
  \Phi^{\boldsymbol{\lambda}}_{\al}=\Phi^{|\boldsymbol{\lambda}|}_\al (x)=\Phi^{\lambda}_\al (x)= {(\lambda )}^{n/4} \Phi_\al({(\lambda ) }^{1/2}x ),
$$
and
$$
\Phi^{\boldsymbol{\lambda}}_{\al,\beta}=\Phi^{|\boldsymbol{\lambda}|}_{\al ,\beta}(x,u)=\Phi^{\la}_{\al,\beta}.
$$
We have the following well known result:
\begin{equation}
    \overline{\Phi^{\la}_{\al ,\beta}}\times_{\la}\overline{\Phi^{\la}_{\mu, \nu}}=
  (2\pi)^{n/2}{(\la)}^{-n}\;\delta_{\nu \al}\;\overline{\Phi^{\la}_{\mu, \beta}}.
\end{equation}
We denote \(\varphi^{n-1}_k\) as the Laguerre function on \(\mathbb{R}^n \times \mathbb{R}^n\). Also, for \(\la >0\), \(\varphi^{n-1}_{k,\la}(x,u)=\varphi^{n-1}_{k,\la}(\sqrt{\la}(x,u)).\)

To make the notations compact, we introduce the following notation, which will be used subsequently while establishing Gutzmer's formula and Beurling's theorem. For \(q,l>0\), \(j,k \in \mathbb{N}\) and \(w=(w',w'')\in \mathbb{C}^{2q}\times \mathbb{C}^{2l}\), let $\varphi^{q,l}_{j,k,\la}(w',w'')=\varphi^{2q-1}_{j,\la}(w')\varphi^{2q-1}_{j,\la}(w'') $.

\subsection{Radon transform}
We recall (\cite{dasgupta-thangavelu-nilmanifold, roncal-thangavelu-holomorphic-extension}) that to every H-type Lie algebra \(\mathfrak{n}=\mathfrak{v}\oplus \mathfrak{z}\) and unit vector \(\omega \in \mathfrak{z}\), one can associate a
 Heisenberg Lie algebra as follows: Let \(k(\omega)\) denote the orthogonal
 complement of \(\omega\) in \(\mathfrak{z}\). Let \(q_{\omega}\) denote the quotient map from \(\mathfrak{n}\) to \(\mathfrak{n}/k(\omega)\). It is known that this quotient algebra can be identified with \(\mathfrak{v}\oplus \mathbb{R}\) and is isomorphic to the Heisenberg Lie algebra \(\mathfrak{h}_{\omega}\). We denote the associated Lie group by \(\mathbb{H}^n_{\omega}\). %This link between H-type Lie algebras and Heisenberg Lie algebras provides insight into the representation theory of H-type groups. As before, we do not consider the one-dimensional representations.  If \(\pi\) is any
 %infinite-dimensional irreducible representation of H-type group \(N\), then its restriction to the center is a unitary character, i.e., \(\exists \mu \in \mathbb{R}^*\) and \(\omega \in \mathbb{S}^{m-1}\), the unit sphere, such that \(\pi(0,t)=e^{i\mu \omega \cdot t}\operatorname{Id}\). Moreover, such a representation factors through
 %a representation of \(\mathbb{H}^n_{\omega}\). By the Stone–von
 %Neumann theorem, all infinite-dimensional irreducible unitary representations of
 %\(N\) are parametrised by \((\mu, \omega)\), \(\mu >0, \omega \in \mathbb{S}^{m-1}\) and we denote such a representation by \(\pi_{\mu,\omega}\). 
For \(\lambda \in \mathbb{R}^*\) and \(\omega \in \mathbb{S}^{m-1}\), it can be seen that the restriction of \(\pi_{\lambda,\omega}\) to \(\mathbb{H}^n_{\omega}\)  is unitarily equivalent to the Schr\"odinger representation \(\pi_\lambda\) of the Heisenberg group. That is, there is an unitary map $\Gamma: L^2(\mathbb R^n)\rightarrow L^2(\mathbb R^n)$ such that 
 \Be 
 \pi_{\lambda,\omega}=\Gamma \,\pi_{\lambda} \,\Gamma^{-1}.
 \Ee
Now, we briefly recall the definition of (partial) Radon transform \cite{dasgupta-thangavelu-nilmanifold, roncal-thangavelu-holomorphic-extension} on \(N\). Given an integrable function \(f\) on \(N\) and \(\omega \in \mathbb{S}^{m-1}\), its (partial) Radon transform is defined as
\begin{equation}\lab{radon transform}
 f_{\omega}(x,u,a)=\int\limits_{t\cdot \omega= a}f(x,u,t)\;d\sigma_1(t),\;\;\; x,u \in \mathbb{R}^n \;\text{and}\; a \in \mathbb{R},
\end{equation}
where \(d\sigma_1\) denotes the Lebesgue measure on the \((m-1)\) dimensional hyperplane \(\{t \cdot \omega=a\}\).
From the definition, it follows that
 $$
 \int\limits_{\mathbb{R}^m}e^{i\lambda \omega \cdot t}f(x,u,t) dt=\int\limits_{\mathbb{R}}e^{i \lambda a}f_{\omega}(x,u,a)da.$$
 That is,  \Be \widetilde{f_{\omega}}(x,u,\lambda)=\widetilde{f}(x,u,\lambda \omega), \Ee where $\widetilde{f}$ and $ \widetilde{f_\omega}$ denote the Fourier transform of $f$ and $f_{\omega}$ in the last variable respectively.
 Also we have, \Be \widetilde{f_{\omega}}(x,u,\lambda)=\widetilde{f_{-\omega}}(x,u,-\lambda). \Ee  
Using the relation between the representation on H-type group and \(\mathbb{H}^n_{\omega}\), it follows that 
 \Be
 \pi_{\lambda,\omega}(f)=\pi_{\lambda}(f_{\omega}).
 \Ee
 Also, it is important to note that the function \(f_{\omega}\) need not belong to \(L^2(\mathbb{H}^n_{\omega})\) but the modified Radon transform \(R_{\omega}f=\left( (-\Delta_t)^{(m-1)/4}f \right)_{\omega}\) does for almost every \(\omega\). In addition, from \cite{dasgupta-thangavelu-nilmanifold}, we have the following relation
\begin{equation}\label{plancherel for radon}
{\lVert f \rVert}_2^2=\int\limits_{\mathbb{S}^{m-1}}\left({\int\limits_{\mathbb{R}^{2n+1}}{|R_{\omega}f(x,u,a)|}^2dxduda}\right) d\omega.
\end{equation}

%\begin{theorem}\label{first Beurling's using radon transform}
%Let \(f\in L^2(N)\) such that \(\mathcal{R}(\Hat{f}(\la,\omega)) \subset \mathcal{D}(\pi_{\la,\omega(z,\gamma))}\) and \(\pi_{\la,\omega}\Hat{f}(\la,\omega)\) is of trace class such that for each \(\omega \in \mathbb{S}^{m-1}\)
%\begin{equation}\label{cond1}
%   \displaystyle\int_{U(n)}{\left \lVert \pi_{\la,\omega}\left({q_{\omega}}^{-1}(\sigma(z,\gamma,0)\right)\Hat{f}(\la,\omega)\right \rVert}^2_{HS}< \infty.
%\end{equation} 
%If for each \(\omega \in \mathbb{S}^{m-1}\) the function satisfies 
 %   \begin{equation}\label{cond2}
  %    \int_{\mathbb{R}^{2n}}\left | f^{\la,\omega}(y,v)  \right | \left ( \sup_{|(y',v')| \leq |(y,v)|} {\left \lVert \pi_{\la,\omega}(i(y',v'))\Hat{f}(\la,\omega)\right \rVert }_1  \right ) dydv < \infty  
   % \end{equation}
   % for almost every \(\la \in \mathbb{R}^*\) then the function \(f=0\).
%\end{theorem}
%\begin{remark}
 %   Note that we have used the quotient map \(q_{\omega}\) because \(\sigma \) acts on elements of the Heisenberg algebra and then using the surjectivity of the quotient map we come back to the H-type algebra. Keeping this in mind, we omit writing \(q_{\omega}\) and simply write \(\pi_{\la,\omega}\left((\sigma(z,\gamma,0)\right) \) instead of \(\pi_{\la,\omega}\left({q_{\omega}}^{-1}(\sigma(z,\gamma,0)\right)\).
%\end{remark}

\section{Proof of Theorem \ref{first Beurling's using radon transform} and Theorem \ref{second Beurling's using radon transform}}\label{sec3}
First, we prove Theorem \ref{first Beurling's using radon transform} below using the Radon transform defined in the section above. 
\begin{remark}
  Note that in the statement of the theorem, we have used the quotient map \(q_{\omega}\) because \(\sigma \) acts on elements of the Heisenberg algebra and then using the surjectivity of the quotient map we come back to the H-type algebra. Keeping this in mind, we omit writing \(q_{\omega}\) and simply write \(\pi_{\la,\omega}\left((\sigma(z,\gamma,0)\right) \) instead of \(\pi_{\la,\omega}\left({q_{\omega}}^{-1}(\sigma(z,\gamma,0)\right)\).  
\end{remark}
\begin{proof}[Proof of Theorem \ref{first Beurling's using radon transform}:]
    Fix an element \(\omega \in \mathbb{S}^{m-1}\). The condition \eqref{cond1} states that 
    \begin{equation*}
   \displaystyle\int_{U(n)}{\left \lVert \pi_{\lambda,\omega}\left((\sigma(z,\gamma)\right)\Hat{f}(\lambda,\omega)\right \rVert}^2_{HS} \,d\sigma< \infty.
\end{equation*} 
We have,
\Be \pi_{\lambda,\omega}(\sigma(z,\gamma))=\Gamma \pi_{\lambda}(\sigma(z,\gamma))\Gamma^{-1}. \Ee
Also, for the Fourier transform of the function \(f\), we have \Be \Hat{f}(\lambda,\omega)=\pi_{\lambda,\omega}(f)=\Gamma\, \pi_{\lambda}(f_\omega)\,\Gamma^{-1}={|\lambda|}^{(1-m)/2}\,\Gamma \,\pi_{\lambda}(R_{\omega}f)\,\Gamma^{-1}.\Ee
Therefore, we get \begin{align*}
 \pi_{\lambda,\omega}(\sigma(z,\gamma))\Hat{f}(\lambda,\omega)= &{|\lambda|}^{(1-m)/2}\,\Gamma\,\pi_{\lambda}(\sigma(z,\gamma))\,\Gamma^{-1}\,\Gamma\,\pi_{\lambda}(R_{\omega}f)\Gamma^{-1}\\
 =&{|\lambda|}^{(1-m)/2}\Gamma\,\pi_{\lambda}(\sigma(z,\gamma))\pi_{\lambda}(R_{\omega}f)\,\Gamma^{-1}.
\end{align*}
Since \(\Gamma\) is unitary, we deduce that \Be
{\left \lVert \pi_{\lambda,\omega}(\sigma(z,\gamma))\Hat{f}(\lambda,\omega) \right \rVert}^2_{HS}={|\lambda|}^{(1-m)}{\left \lVert \pi_{\lambda}(\sigma(z,\gamma))\pi_{\lambda}(R_{\omega}f) \right \rVert}^2_{HS}. \Ee
Observe that we are taking the modified Radon transform \(R_{\omega}f\) in the picture because \(f_{\omega}\) need not be in \(L^2(\mathbb{H}^n_{\omega})\).
Thus the condition \eqref{cond1} becomes, 
\begin{equation}\label{cond3}
   \displaystyle\int_{U(n)}{\left \lVert \pi_{\lambda}\left((\sigma(z,\gamma)\right)\widehat{R_{\omega}f}(\lambda)\right \rVert}^2_{HS} \,d\sigma < \infty.
\end{equation}
The condition \eqref{cond2} states that 
\begin{equation*}
    \int_{\mathbb{R}^{2n}}\left | f^{\lambda\omega}(y,v)  \right | \left ( \sup_{|(y',v')| \leq |(y,v)|} {\left \lVert \pi_{\lambda,\omega}(i(y',v'))\Hat{f}(\lambda,\omega)\right \rVert }_1  \right ) dydv < \infty.
\end{equation*}
Using the relation between \(\pi_{\lambda,\omega}\) and \(\pi_{\lambda}\), we obtain the following
\begin{align*}
\sup_{|(y',v')| \leq |(y,v)|} {\left \lVert \pi_{\lambda,\omega}(i(y',v'))\Hat{f}(\lambda,\omega)\right \rVert }_1&=\sup_{|(y',v')| \leq |(y,v)|} {\left \lVert \pi_{\lambda,\omega}(i(y',v'))\pi_{\lambda,\omega}(f)\right \rVert }_1\\
&=\sup_{|(y',v')| \leq |(y,v)|} {\left \lVert \Gamma \pi_{\lambda}(i(y',v'))\Gamma^{-1} \left( {|\lambda|}^{(1-m)/2}\Gamma \pi_{\lambda}(R_{\omega}f) \Gamma^{-1} \right)\right \rVert  }_1\\
&=\sup_{|(y',v')| \leq |(y,v)|} {|\lambda|}^{(1-m)/2}{\left \lVert \Gamma \pi_{\lambda}(i(y',v'))  \pi_{\lambda}(R_{\omega}f) \Gamma^{-1}\right \rVert }_1\\
&=\sup_{|(y',v')| \leq |(y,v)|} {|\lambda|}^{(1-m)/2}{\left \lVert \pi_{\lambda}(i(y',v'))  \pi_{\lambda}(R_{\omega}f) \right \rVert }_1.
\end{align*}
Moreover, from the definition of \(f^{\lambda\omega}\), we deduce that
$$
f^{\lambda\omega}(y,v)=\int\limits_{\mathbb{R}^m}f(y,v,t)e^{i\lambda \omega \cdot t}dt=\int\limits_{\mathbb{R}}f_{\omega}(a)e^{i\lambda a}da=(f_{\omega})^{\lambda}(y,v).
$$
Since \((R_{\omega}f)(y,v,a)=D_{a}^{(m-1)/2}f_{\omega}(y,v,a)\), where \(D_{a}^{(m-1)/2}\) is the fractional derivative of order \((m-1)/2\), we get, \Be (R_{\omega}f)^{\lambda}(y,v)={|\lambda|}^{(m-1)/2}(f_{\omega})^{\lambda}(y,v)={|\lambda|}^{(m-1)/2}f^{\lambda\omega}(y,v).\Ee
Thus the condition \eqref{cond2} becomes
\begin{equation*}
   {|\lambda|}^{(1-m)} \int_{\mathbb{R}^{2n}}\left | (R_{\omega}f)^{\lambda}(y,v)  \right | \left ( \sup_{|(y',v')| \leq |(y,v)|} {\left \lVert \pi_{\lambda}(i(y',v'))\widehat{(R_{\omega}f)}(\lambda)\right \rVert }_1  \right ) dydv < \infty.
\end{equation*}
That is, for almost every \(\lambda \in \mathbb{R}^*\),
\begin{equation}\label{cond4}
   \int_{\mathbb{R}^{2n}}\left | (R_{\omega}f)^{\lambda}(y,v)  \right | \left ( \sup_{|(y',v')| \leq |(y,v)|} {\left \lVert \pi_{\lambda}(i(y',v')) \widehat{(R_{\omega}f)}(\lambda)\right \rVert }_1  \right ) dydv < \infty.
\end{equation}
Now combining Eq. \eqref{cond3} and Eq. \eqref{cond4}, we have that for the function \(R_{\omega}f \in L^2(\mathbb{H}^n_{\omega})\),
$$
\displaystyle\int_{U(n)}{\left \lVert \pi_{\lambda}\left((\sigma(z,\gamma)\right)\widehat{R_{\omega}f}(\lambda)\right \rVert}^2_{HS} \,d\sigma< \infty$$
and
$$
\int_{\mathbb{R}^{2n}}\left | (R_{\omega}f)^{\lambda}(y,v)  \right | \left ( \sup_{|(y',v')| \leq |(y,v)|} {\left \lVert \pi_{\lambda}(i(y',v'))\widehat{(R_{\omega}f)}(\lambda)\right \rVert }_1  \right ) dydv < \infty.
$$
Now by applying Theorem \ref{beurling on heisenberg} on the function \(R_{\omega}f\), we conclude that \(R_{\omega}f=0\).
This is true for all \(\omega \in \mathbb{S}^{m-1}\). Lastly, using Eq. \eqref{plancherel for radon}, we conclude that the function \(f=0\).
\end{proof}

Next we shall prove Theorem \ref{second Beurling's using radon transform}.
\begin{proof}[Proof of Theorem \ref{second Beurling's using radon transform}:] The proof follows along the proof of \cite[Theorem 2.3]{Beurling-rudra-parui}.
  Fix \(\omega \in \mathbb{S}^{m-1}\). We can assume that \({\left \lVert \Hat{f}(\lambda,\omega)\right \rVert}_{HS}|\lambda|^n\) is different from \(0\) on a set of positive measure. Thus we have 
  \begin{equation}\label{second Beurling's eq2}
      \int\limits_{\mathbb{C}^{n}\times \mathbb{R}^m}\frac{\left \lvert (-\Delta_t)^{(m-1)/4}f(z,t)\right \rvert e^{\lvert t \rvert \lvert  y_0|}}{{(1+|z|)}^{M_1}{(1+|y_0|+|t|)}^{M_2}}dz dt<\infty,
  \end{equation}
  for \(y_0\) in a set of positive measure in \(\mathbb{R}\).
In Eq. \eqref{second Beurling's eq1}, using the relation between \(\pi_{\lambda,\omega}f\) and \(\pi_{\lambda}(R_{\omega}f)\), we get that
\begin{equation}\label{second Beurling's eq3}
\int\limits_{\mathbb{R}^{2n}\times \mathbb{R}^m}{\int\limits_{\mathbb{R}}\frac{\left \lvert (-\Delta_t)^{(m-1)/4}f(z,t)\right \rvert{\lVert \pi_{\lambda}(R_{\omega}f)\rVert}_{HS}e^{\lvert t \rvert \lvert  \lambda|}{\lvert \lambda \rvert}^{n}}{{(1+|z|)}^{M_1}{(1+|\lambda|+|t|)}^{M_2}}d\lambda dz dt}< \infty. 
\end{equation}
We claim that for \(\omega_1 \in \mathbb{S}^{m-1}\), the following holds
\begin{equation}\label{claim integral}
I=\int\limits_{\mathbb{R}^{2n}\times \mathbb{R}}{\int\limits_{\mathbb{R}}\frac{\left \lvert R_{\omega_1}f(z,a)\right \rvert{\lVert \pi_{\lambda}(R_{\omega}f)\rVert}_{HS}e^{\lvert a \rvert \lvert  \lambda|}{\lvert \lambda \rvert}^{n}}{{(1+|z|)}^{M_1}{(1+|\lambda|+|a|)}^{M_2}}d\lambda dz da}< \infty. 
\end{equation}
To prove the above claim, we break the integral $I$ in three parts. Let,
\begin{align*}
    &I_1=\int\limits_{\mathbb{R}^{2n}}\int\limits_{\mathbb{R}}{\int\limits_{|\lambda|>L}\frac{\left \lvert R_{\omega_1}f(z,a)\right \rvert{\lVert \pi_{\lambda}(R_{\omega}f)\rVert}_{HS}e^{\lvert a \rvert \lvert  \lambda|}{\lvert \lambda \rvert}^{n}}{{(1+|z|)}^{M_1}{(1+|\lambda|+|a|)}^{M_2}}d\lambda dz da},\;\; \text{where }L+L^2>M_2,\\
    &I_2=\int\limits_{\mathbb{R}^{2n}}\int\limits_{|a|>M}{\int\limits_{|\lambda|\leq L}\frac{\left \lvert R_{\omega_1}f(z,a)\right \rvert{\lVert \pi_{\lambda}(R_{\omega}f)\rVert}_{HS}e^{\lvert a \rvert \lvert  \lambda|}{\lvert \lambda \rvert}^{n}}{{(1+|z|)}^{M_1}{(1+|\lambda|+|a|)}^{M_2}}d\lambda dz da}\;\;\text{where }M=2(L+1),\\
    &I_3=\int\limits_{\mathbb{R}^{2n}}\int\limits_{|a|\leq M}{\int\limits_{|\lambda|\leq L}\frac{\left \lvert R_{\omega_1}f(z,a)\right \rvert{\lVert \pi_{\lambda}(R_{\omega}f)\rVert}_{HS}e^{\lvert a \rvert \lvert  \lambda|}{\lvert \lambda \rvert}^{n}}{{(1+|z|)}^{M_1}{(1+|\lambda|+|a|)}^{M_2}}d\lambda dz da}.
\end{align*}
Following the lines of proof of \cite[Theorem 2.3]{Beurling-rudra-parui}, consider the function 
$$
F(z)=\frac{e^{\alpha z}}{{(1+\al+z)}^{M_2}}\;\; \text{for }\al>0\text{ and }\al+\al^2>M_2.
$$
Then if \(z_1\geq z_2\), \(F(z_1)\geq F(z_2)\). So, for \(t\in \mathbb{R}^m\), taking \(z_1=|t|\), \(z_2=|\langle t, \omega_1 \rangle|\) and \(\al=|\lambda|\geq L\), we get \Be \frac{e^{|\lambda||t|}}{{(1+|\lambda|+|t|)}^{M_2}}\geq \frac{e^{|\lambda||\langle t, \omega_1 \rangle|}}{{(1+|\lambda|+|\langle t, \omega_1 \rangle|)}^{M_2}}.\Ee  Applying this in Eq. \eqref{second Beurling's eq3}, we obtain 
\begin{equation*}
\int\limits_{\mathbb{R}^{2n}}\int\limits_{\mathbb{R}}{\int\limits_{t \cdot \omega_1=a}{\int\limits_\mathbb{R}}\frac{\left \lvert (-\Delta_t)^{(m-1)/4}f(z,t)\right \rvert{\lVert \pi_{\lambda}(R_{\omega}f)\rVert}_{HS}e^{|\lambda||\langle t, \omega_1 \rangle|}{\lvert \lambda \rvert}^{n}}{{(1+|z|)}^{M_1}{(1+|\lambda|+|\langle t, \omega_1 \rangle|)}^{M_2}}d\lambda d\sigma_1(w_1) da dz}<\infty,
\end{equation*}
where \(d\sigma_1\) denotes the Lebesgue measure on the hyperplane \(\{t: t\cdot \omega_1=a\}\). From this, we get
\begin{equation*}
\int\limits_{\mathbb{R}^{2n}}\int\limits_{\mathbb{R}}{\int\limits_{|\lambda|>L}\frac{\left \lvert R_{\omega_1}f(z,a)\right \rvert{\lVert \pi_{\lambda}(R_{\omega}f)\rVert}_{HS}e^{\lvert \lambda \rvert \lvert  a|}{\lvert \lambda \rvert}^{n}}{{(1+|z|)}^{M_1}{(1+|\lambda|+|a|)}^{M_2}}d\lambda  da dz}<\infty.
\end{equation*}
Thus \(I_1<\infty\).

Now we consider $I_2$.
We recall that from the Plancherel theorem of the Heisenberg group, for a square-integrable function \(g\) on  \(\mathbb{H}^n\), we have \Be
{\lVert f \rVert}^2_{2}={(2\pi)}^{(-n-1)} \displaystyle\int{{\lVert \Hat{f}(\lambda)\rVert}^2_{HS} |\lambda|^nd\lambda}. \Ee  Thus, we have, 
\begin{equation}\label{second Beurling's eq4}
I_2\leq c \int\limits_{\mathbb{R}^{2n}}\int\limits_{|a|>M}\frac{|R_{\omega_1}f(z,a)|e^{|a|L}}{{(1+|z|)}^{M_1}{(1+|a|)}^{M_2}}da dz.
\end{equation}
Note that in Eq. \eqref{second Beurling's eq2}, the point \(y_0\) can be chosen such that \(|y_0|>2L\), so that,
\begin{equation*}
      \int\limits_{\mathbb{R}^{2n}\times \mathbb{R}^m}\frac{\left \lvert (-\Delta_t)^{(m-1)/4}f(z,t)\right \rvert e^{\lvert t \rvert \lvert  y_0|}}{{(1+|z|)}^{M_1}{(1+|y_0|+|t|)}^{M_2}}dz dt<\infty.
  \end{equation*}
Now similar to above, using the function \(F(z)\), with \(\al=|y_0|>2L\), \(z_1=|t|\) and \(z_2=|\langle t, \omega_1 \rangle|\) we get
\begin{equation*}  \int\limits_{\mathbb{R}^{2n}}\int\limits_{\mathbb{R}}\int\limits_{t\cdot \omega_1=a}\frac{\left \lvert (-\Delta_t)^{(m-1)/4}f(z,t)\right \rvert e^{\lvert  y_0|\lvert \langle t, \omega_1 \rangle \rvert }}{{(1+|z|)}^{M_1}{(1+|y_0|+|\langle t, \omega_1 \rangle|)}^{M_2}}d\sigma_1(w_1) da dz <\infty.
\end{equation*}
Therefore, 
\begin{equation*}  \int\limits_{\mathbb{R}^{2n}}\int\limits_{|a|>M}\frac{\left \lvert R_{\omega_1}f(z,a)\right \rvert e^{\lvert  y_0|\lvert a \rvert }}{{(1+|z|)}^{M_1}{(1+|y_0|+|a|)}^{M_2}}da dz<\infty.
\end{equation*}
Again, take \(\al=|a|>M\) and \(z_1=|y_0|,z_2=2L\) to obtain
\begin{equation}\label{second Beurling's eq5}  \int\limits_{\mathbb{R}^{2n}}\int\limits_{|a|>M}\frac{\left \lvert R_{\omega_1}f(z,a)\right \rvert e^{\lvert a \rvert 2L}}{{(1+|z|)}^{M_1}{(1+|a|+2L)}^{M_2}}da dz<\infty.
\end{equation}
Furthermore, it is easy to check that \(\frac{e^{|a|2L}}{{(1+|a|+2L)}^{M_2}} \geq \frac{e^{L|a|}}{{(1+|a|)}^{M_2}}\) if \(|a|>M\). Consequently using this in Eq. \eqref{second Beurling's eq5}, we get
\begin{equation*} \int\limits_{\mathbb{R}^{2n}}\int\limits_{|a|>M}\frac{\left \lvert R_{\omega_1}f(z,a)\right \rvert e^{\lvert a \rvert L}}{{(1+|z|)}^{M_1}{(1+|a|)}^{M_2}}da dz<\infty.
\end{equation*}
Subsequently, with the help of the above estimate and Eq. \eqref{second Beurling's eq4}, we conclude that \(I_2< \infty\).

Now we shall prove that $I_3<\infty$.
By definition of integral \(I_3\), we have
\begin{align*}
    I_3&=\int\limits_{\mathbb{C}^{n}}\int\limits_{|a|\leq M}{\int\limits_{|\lambda|\leq L}\frac{\left \lvert R_{\omega_1}f(z,a)\right \rvert{\lVert \pi_{\lambda}(R_{\omega}f)\rVert}_{HS}e^{\lvert a \rvert \lvert  \lambda|}{\lvert \lambda \rvert}^{n}}{{(1+|z|)}^{M_1}{(1+|\lambda|+|a|)}^{M_2}}d\lambda dz da}\\
    & \leq c \int\limits_{\mathbb{C}^{n}}\int\limits_{|a|\leq M}{\frac{\left \lvert R_{\omega_1}f(z,a)\right \rvert}{{(1+|z|)}^{M_1}} da dz }\\
    & \leq c \left ( \int\limits_{\mathbb{C}^{n}}\int\limits_{|a|\leq M}{\left \lvert R_{\omega_1}f(z,a)\right \rvert}^2 da dz  \right)^{1/2}\; \left ( \int\limits_{\mathbb{C}^{n}}\int\limits_{|a|\leq M}{\left \lvert \frac{1}{{(1+|z|)}^{M_1}}\right \rvert}^2 da dz  \right)^{1/2} \\
    %&  \leq c_1 \int\limits_{\mathbb{C}^{n}}\int\limits_{\mathbb{R}}{\left \lvert R_{\omega_1}f(z,\xi)\right \rvert d\xi dz }\\ 
   % &\leq c\int\limits_{\mathbb{C}^{n}}\int\limits_{\mathbb{R}}\int\limits_{t \cdot \omega_1=\xi}{\left \lvert (-\Delta_t)^{(m-1)/4}f(z,t) \right \rvert d\sigma_1 d\xi dz }\\
   % & = c_1 {\lVert (-\Delta_t)^{(m-1)/4}f \rVert}_2 <\infty.
   & < \infty.
\end{align*}
This proves our claim in \eqref{claim integral}. 
Now we apply Theorem \ref{Beurling-Rudra-Parui}, on the integral \(I\) and obtain our desired result.
\end{proof}
\begin{remark}
    Let \(f\in L^1\cap L^2(N)\) and \((-\Delta_t)^{(m-1)/4}f\in L^1\cap L^2(N)\). If Eq. \eqref{second Beurling's eq1} holds for any  \(M_1>0, M_2\geq 0\) then the conclusion of Theorem \ref{second Beurling's using radon transform} holds true. In the proof above we can see that the given condition implies $I_3$ is finite and rest other part of the proof is same.
\end{remark}
\section{Gutzmer's Formula}\label{sec4}
 We shall now concentrate towards Theorem \ref{Beurling's theorem}. To prove this theorem (following method of Thangavelu) we aim to prove (in this section) Gutzmer's formula.  We first define the H-type motion group in this section, utilizing the automorphism group of \(N\).

Let \(A(N)\) denote the group of automorphisms of \(N\) that act as orthogonal transformations on the Lie algebra \(\mathfrak{n}\). The group \(A(N)\) can be seen as a product of two groups described below. 
Let \(U\) denote the subgroup of automorphisms that act trivially on the center, i.e., 
$$
U=\{\theta \in A(N): \theta(z)=z,\; \forall z \in \mathfrak{z}\}.
$$
For each unit vector $z \in \mathfrak{z}$, the map $J_z: \mathfrak{z}^\perp \rightarrow \mathfrak{z}^\perp$ can be extended to an automorphism of $\mathfrak{n}$ by defining it as the following operator
    $$
     J_{z}(w)= \begin{cases}
     z &\textrm{if}\;\;w=z\\
     -w&\textrm{if } w\in \mathfrak{z} \textrm{ and }\left \langle z,w\right \rangle =0.
     \end{cases}
$$
Let \(\operatorname{Pin}(m)\) denote the subgroup of \(A(N)\) generated by \(\{J_z: z \in \mathfrak{z}\}\). Then the subgroups \(U\) and \(\operatorname{Pin}(m)\) commute and $A(N)=U \cdot \operatorname{Pin}(m)$ unless $m \equiv 1(\bmod 4)$ and in that case $A(N) /(U \cdot \operatorname{Pin}(m))$ has two elements. For more details, one can refer to \cite{kaplanricci,riccicomm}.

Recall that, a pair \((K,N)\) where \(K\) is a subgroup of group of automorphisms of a Lie group \(N\), is said to be a Gelfand pair if the algebra $L_K^1(N)$ of $K$-invariant integrable functions on $N$ forms a commutative algebra under usual convolution. For H-type groups, we have the following classification of Gelfand pairs \cite{riccicomm, yasser2024gelfand}.
\begin{theorem}
    The groups $N$ of H-type for which $L_{A(N)}^1(N)$ is commutative, that is $(A(N), N)$ is a Gelfand pair, are those for which
$$
\operatorname{dim}(\mathfrak{z})=m=\left\{\begin{array}{l}
1,2 \text { or } 3 \\
5,6 \text { or } 7 \text { and } \mathfrak{z}^\perp \text { is irreducible } \\
7, \mathfrak{z}^\perp \text { is isotypic and } \operatorname{dim}(\mathfrak{z}^\perp)=16 .
\end{array}\right.
$$
\end{theorem}
\begin{corollary}
     Let $U$ be the subgroup of $A(N)$ consisting of automorphisms that act trivially on the center $\mathfrak{z}$. Then $(U, N)$ is a Gelfand pair if and only if $\operatorname{dim}{\mathfrak{z}}=1,2$, or 3 .
\end{corollary}
We refer the reader to \cite{riccicomm} for more details about the automorphism group of \(N\). In this article, we will restrict ourselves to the case when \(m=\operatorname{dim}(\mathfrak{z})=1,2 \text{ or }3\). 
\begin{itemize}
    \item \(m=1\), \(N\) is the Heisenberg group, \(\mathbb{H}^n\) and \(U\) is the group \(U(n)\).
    \item \(m=2\), \(\mathfrak{v}=\mathfrak{z}^\perp \cong \mathbbmss{H}^q  \cong \mathbb{C}^{2q}\) , where \(\mathbbmss{H}\) refers to the space of quaternions, \(U\) is the compact sympletic group \(Sp(q)\). Thus in this case, \(N= \mathbb{C}^{n} \times \mathbb{R}^2\), with $n=2q$.
    \item \(m=3\), \(\mathfrak{v}=\mathfrak{z}^\perp \cong \mathbbmss{H}^q \times \mathbbmss{H}^l \cong \mathbb{C}^{2q} \times \mathbb{C}^{2l}\), \(k,l \geq 0\) and the subgroup \(U \cong Sp(q) \times Sp(l)\). Thus in this case, \(N=\mathbb{C}^{n} \times \mathbb{R}^3\), with $n=2q+2l$.
\end{itemize}
Note that here the compact symplectic group refers to the group \(Sp(q)=Sp(2q,\mathbb{C}) \cap U(2q)\), where \(U(2q)\) is the set of \(2q \times 2q\) unitary matrices and
$$
Sp(2q,\mathbb{C})=\{ X \in M_{2q \times 2q}(\mathbb{C}): X^t\Omega X=\Omega\}, 
$$
where \( \Omega =
\begin{bmatrix}
0 & I_q \\
-I_q & 0
\end{bmatrix}
\).

 From now on we take an H-type group $N$ with the dimension of the center $m=2, 3$ and a subgroup 
\(U\) of \(A(N)\) that acts trivially on the center. That is, 
 the action of \(U\) on \(N\) is as follows:
$$
\sigma(z,t)=(\sigma z,t), \;\; \forall \sigma \in U, (z,t) \in N.
$$

 %The Heisenberg motion group and its representation theory are briefly discussed in \cite{thangaveluheisenberg}. Drawing inspiration from this foundational work, we define the H-type motion group in this section, utilizing the automorphism group of \(N\). Furthermore, we present a family of class-one representations of this group. For this discussion, \(N\) is understood to denote the H-type group throughout this section.
 
 We define the {\em H-type motion group} as the semi-direct product of \(U\) and \(N\) \((U \ltimes N)\). %and denote it by \(N_m\).
   The group law is given by
   \Be
   (\sigma,z,t)(\tau,w,s)=\left( \sigma \tau,z+\sigma w, t+s+\frac{1}{2}[z,\sigma w]\right).
   \Ee
   Furthermore, the H-type motion group acts on \(N\) as
   \Be
   (\sigma,z,t)(w,s)=\left(z+\sigma w,t+s+\frac{1}{2}[z,\sigma w]\right).
   \Ee
To investigate the connection between the Fourier transform on the H-type group and the  H-type motion group, we define a family of class one representations of the H-type motion group.

We will associate H-type group \(N\) with \(\mathbb{R}^n \times \mathbb{R}^n \times \mathbb{R}^m\) and write \((x,y,t)\) for \((x+iy,t)\). Moreover, each element \(\sigma\) in \(U\) can be written as \(A+iB\) where \(A\) and \(B\) are real and imaginary parts of \(\sigma\). Then the action of \(U\) on \(\mathbb{R}^n \times \mathbb{R}^n \times \mathbb{R}^m\) can be written as 
$$
\sigma(x,u,t)=(Ax-By,Ay+Bx,t).
$$
We extend this to an action of \(U\) on \(\mathbb{C}^n \times \mathbb{C}^n \times \mathbb{C}^m\) as
\begin{equation}\label{ext-action-of-U}
\sigma(z,w,\zeta)=(Az-Bw,Aw+Bz,\zeta), \;\; \forall z,w \in \mathbb{C}^n \text{ and } \zeta \in \mathbb{C}^m.
\end{equation}
We recall that, the H-type motion group acts on \(\mathbb{R}^n \times \mathbb{R}^n \times \mathbb{R}^m\) as
$$
(\sigma,z,t)(w,s)=(z+\sigma w,t+s+\frac{1}{2}[z,\sigma w]).
$$
Using \eqref{ext-action-of-U}, this action can be extended to \(\mathbb{C}^n \times \mathbb{C}^n \times \mathbb{C}^m\) by
$$
(\sigma,x,u,t)(z,w,\zeta)=(x,u,t)(Az-Bw,Bz+Aw,\zeta) 
$$
where \(\sigma=A+iB\).
Furthermore, this action extends to functions defined on \(\mathbb{C}^n \times \mathbb{C}^n \times \mathbb{C}^m\) as
$$
\rho(g)f(z,w,\zeta)=f(g^{-1}(z,w,\zeta)), \; \text{where } g \in U \ltimes N.
$$

%To prove Gutzmer's formula for the H-type group $N$, we need to concentrate only on the case when $(U, N)$ is a Gelfand pair. Therefore we will concentrate on the cases when $m=1, 2, 3$. But the case when $m=1$, corresponds to the Heisenberg group. Therefore 

We now consider the two cases $m=2$ and $m=3$ separately. First we consider the case when $m=3$. 

 \subsection{The case when $m=3$:}
  We recall that in this case, \(N=\mathbb{C}^{2q} \times \mathbb{C}^{2l} \times \mathbb{R}^3\) and \(U= Sp(q) \times Sp(l) \). The group \(U\) acts transitively on the product of spheres and we know H-type motion group is \(U \ltimes N\).%Also, \(N_3\) is the semidirect product of \(Sp(q) \times Sp(l) \) and \(N\).
  
 For \(j,k \in \mathbb{N},\; \boldsymbol{\lambda}=\la \omega \in \mathbb{R}^3\setminus \{0\}\), let \(\mathcal{H}^{\boldsymbol{\lambda}}_{j,k}\) be the Hilbert space with orthonormal basis 
 \Be
 E^{\boldsymbol{\lambda}}_{\al,\beta}(z,t)=e^{-i \boldsymbol{\lambda} \cdot t}\Phi^{\la}_{\al',\beta'}(z')\Phi^{\la}_{\al'',\beta''}(z''), \;\;t \in \mathbb{R}^3 
 \Ee
 where \(\la=|\boldsymbol{\la}|\), \(\al=(\al',\al'') \in \mathbb{N}^{2q} \times \mathbb{N}^{2l},\beta=(\beta', \beta'')\in \mathbb{N}^{2q} \times \mathbb{N}^{2l} \) with \(|\al'|=j\), \(|\al''|=k\) and \(z=(z',z'') \in \mathbb{C}^{2q}\times \mathbb{C}^{2l}\). 
Note that here, \(\Phi^{\la}_{\al',\beta'}(z')\) and \(\Phi^{\la}_{\al'',\beta''}(z'')\) refer to the special Hermite function on \(\mathbb{C}^{2q}\) and \(\mathbb{C}^{2l}\) respectively.

 We equip the Hilbert space \(\mathcal{H}^{\boldsymbol{\lambda}}_{j,k}\) with the following inner product structure
 $$
\left \langle F,G\right \rangle= {( \la )}^{2q+2l} \int\limits_{\mathbb{C}^{2q} \times \mathbb{C}^{2l}}{F(z',z'',0)\;\overline{G(z',z'',0)} \;dz'dz''}.
$$
We recall that \(n=2q+2l\). The elements of the space \(\mathcal{H}^{\boldsymbol{\lambda}}_{j,k}\) can be characterised by the following eigenvalue equations $$
\mathcal{L}(f)=(2j+2k+n)\lvert \boldsymbol{\lambda} \rvert f=(2j+2k+n)\la f, \;\;\;i\frac{\partial}{\partial t_j}=\la_j f, \;\; 1 \leq j \leq 3.
$$
Moreover, if 
\Be
e^{-i \boldsymbol{\lambda} \cdot t}f(z',z'') \in \mathcal{H}^{\boldsymbol{\lambda}}_{j,k},
\Ee
then for each \(z'' \in \mathbb{C}^{2l}\), \(f(z',z'')\) is an eigenfunction of \(\mathcal{L}_{\la}\) on \(\mathbb{C}^{2q}\). Similarly, for each \(z' \in \mathbb{C}^{2q}\), \(f(z',z'')\) is an eigenfunction of \(\mathcal{L}_{\la}\) on \(\mathbb{C}^{2l}\).

 We define a family \(\rho^{\boldsymbol{\lambda}}_{j,k} \), of class \(1\) representations  of the H-type motion group  on the space \(\mathcal{H}^{\la}_{j,k}\) as follows:
\Be
\rho^{\boldsymbol{\lambda}}_{j,k}\left((\sigma',\sigma''),(z',z''),t\right)\varphi\left((w',w''),s\right)=\varphi \left( \left((\sigma',\sigma''),(z',z''),t\right)^{-1}\left((w',w''),s\right)\right),
\Ee
where $\varphi \in \mathcal{H}^{\boldsymbol{\lambda}}_{j,k}.$
We now prove the irreducibility of the representation defined above. %\textcolor{red}{We note that in this case, the group \(U=Sp(q) \times Sp(l)\) acts transitively on the product of spheres centered at the origin.}
\begin{theorem}
    The representation \(\rho^{\boldsymbol{\lambda}}_{j,k} \) is an irreducible unitary representation of the H-type motion group on \(\mathcal{H}^{\boldsymbol{\lambda}}_{j,k}\).
\end{theorem} 
\begin{proof}
%For this particular proof, we identify the Lie algebra of \(N\) with \(\mathbb{C}^n \times \mathbb{R}^3\) with the specific basis (as mentioned in Section \ref{sec2}) corresponding to \(-\la\), where \(\la \in \mathbb{R}^3\). 
Let \(M\) be a closed invariant subspace of \(\mathcal{H}^{\boldsymbol{\lambda}}_{j,k}\). We will show that \(M=\mathcal{H}^{\boldsymbol{\lambda}}_{j,k}\). Let 
    $$
    M_1=\{f(w)=f(w',w''): e^{-i \boldsymbol{\lambda} \cdot s} \overline{f(w',w'')}\in M\}.
    $$
   % Note that \(M_1\) is a closed subspace of \(L^2(\mathbb{C}^n)\). Also, for any \(f \in M_1\), \(\overline{f}= \sum\limits_{| \al |=k}{\sum\limits_{\beta}{c_{\al \beta}\Phi^{|\la|}_{\al ,\beta}}}\). So we obtain that \(f \times_{-\lvert \la \rvert} \varphi^\la_k=(2\pi)^n {\lvert \la \rvert }^{-n}f\).
 Now, \(M\) is an invariant subspace of \(\mathcal{H}^{\boldsymbol{\lambda}}_{j,k}\). So,
    $$
    \rho^{\boldsymbol{\lambda}}_{j,k}(\sigma,z,t)M \subseteq M, \;\; \forall \;\sigma=(\sigma',\sigma'') \in U,\;z=(z',z'') \in \mathbb{C}^{2q} \times \mathbb{C}^{2l},\;t \in \mathbb{R}^3 .
    $$
    In particular, \(\rho^{\boldsymbol{\lambda}}_{j,k}((\sigma,z,t)^{-1})
    F(w,s) \subseteq M\) for every function \(F(w,s) \in M\). If \(f \in M_1\), that is, 
    \(F(w,s)=e^{-i \boldsymbol{\lambda} \cdot s}\overline{f(w)} \in M\). Then, the condition above can be written as 
    $$
    e^{-i \boldsymbol{\lambda}\cdot \left ( t+s+\frac{1}{2}[z,\sigma w]\right )}\overline{f(z+\sigma w)} \in M.
    $$
   In particular, this implies that
    $$
    e^{\frac{i}{2}\boldsymbol{\lambda}\cdot [z ,\sigma w]}f(z+\sigma w)\in M_1 
    $$
    and
    $$
    f(-w) \in M_1,
    $$
   for all \( z \in \mathbb{C}^n\) and \(\sigma \in U\).
  Next, we consider the following function of the subspace \(M_1\), given by 
  $$
  f_z(w)=\int\limits_{U}f(z+\sigma w)e^{i \frac{\boldsymbol{\lambda}}{2} \cdot [z,\sigma w]}d \sigma.
  $$
  Observe that there exists some \(z\) for which that function \(f_z(w)\) is not identically zero. If \(f_z(w)\) is identically for all \(z,w\), then as \(w=(w',w'') \rightarrow 0\),
  $$
  0=\int\limits_{U}f(z+\sigma w)e^{i\frac{\boldsymbol{\lambda}}{2} \cdot [z,\sigma w]}d \sigma \rightarrow\int\limits_{U} f(z) d\sigma= f(z).
  $$
  This implies that the function \(f\) is identically 0. This is a contradiction as the set \(M_1\) contains non-trivial functions.
  
  By the transitive action of the group \(U\) on product of spheres, the function \(f_z(w)=f_z(w',w'')\) is a polyradial eigenfunction, that is, it is radial in \(w'\) and \(w''\). By Lemma \ref{twisted laplacian lemma}, the action of sublaplacian on function \(e^{-i \boldsymbol{\lambda} \cdot t}f_z(w',w''))\) is given by the sum of the operators \(\mathcal{L}_{\lambda}\) acting on \(\mathbb{C}^{2q}\) and \(\mathcal{L}_{\lambda}\) acting on \(\mathbb{C}^{2l}\). 
  More explicitly, if we consider the function of the form \(g(w',w'',t)=g(r,s,t)\), that is radial in \(z'\) and \(z''\) then the action of the sublaplacian \(\mathcal{L}\) reduces to the operator \(P\) acting on \(g\) given by 
  $$
  P=-\left ( \frac{\partial^2}{\partial r^2}+ \frac{\partial^2}{\partial s^2}\right)-\left( \frac{4q-1}{r}\frac{\partial }{\partial r}+\frac{4l-1}{s}\frac{\partial}{\partial s}\right)+\frac{1}{4}\left( r^2+s^2\right)T\;+P',
  $$
  where
  $$
P'=-\sum_{k=1}^3\left(\sum_{i=1}^{2q}\sum_{j=1}^{2l}w'_i{w''}_jU^{k}_{i,2q+j}\left ( \frac{\partial^2}{\partial s \partial t_k}-\frac{\partial^2}{\partial r \partial t_k}\right )\right ),
  $$
and \(U^k\) are matrices as in equation \eqref{group law}.\\
  Now, \(f_z\) is radial in \(w'\) and \(w''\) respectively and it is well known that the radial eigenfunctions of \(\mathcal{L}_{\lambda}\) are given by Laguerre function. Thus, \Be 
  f_z(w)=f_z(w',w'')=c\, \varphi^{2q-1}_{j,\la}(w')\varphi^{2l-1}_{k,\la}(w'')= \varphi^{q,l}_{j,k,\la}(w',w''),\Ee  for some constant \(c\), where \(\varphi^{2q-1}_{j,\la}, \varphi^{2l-1}_{k,\la}\) are Laguerre functions on \(\mathbb{C}^{2q}\) and \(\mathbb{C}^{2l}\) respectively.
  % Note that by the transitive action of the group \(U\) on the product of spheres centered at the origin, we obtain that the function \(f_z(w)\) is a polyradial eigenfunction of the twisted Laplacian \(\mathcal{L}_{\la}\) on \(\mathbb{C}^n\). Also, for each \(w'' \in \mathbb{C}^{2l}, f_{z}(w',w'')\) is a radial eigenfunction of \(\mathcal{L}_{\la}\) on \(\mathbb{C}^{2q}\) and for each \(w' \in \mathbb{C}^{2q}, f_{z}(w',w'')\) is a radial eigenfunction of \(\mathcal{L}_{\la} \) on \(\mathbb{C}^{2l}\). Thus, \Be 
%  f_z(w)=f_z(w',w'')=c\, \varphi^{2q-1}_{j,\la}(w')\varphi^{2l-1}_{k,\la}(w'')= \varphi^{q,l}_{j,k,\la}(w',w''),\Ee  for some constant \(c\), where \(\varphi^{2q-1}_{j,\la}, \varphi^{2l-1}_{k,\la}\) are Laguerre functions on \(\mathbb{C}^{2q}\) and \(\mathbb{C}^{2l}\) respectively. %\(\mathbb{R}^q \times \mathbb{R}^q\) and \(\mathbb{R}^l \times \mathbb{R}^l\) respectively. 
  So the function \Be \psi(w',w'')=\varphi^{q,l}_{j,k,\la}(w',w'') \in M_1. \Ee
  % Next, we claim that \(\varphi^{|\la|}_k \in M_1\). If not then, \(\exists\; h \in L^2(\mathbb{C}^n)\) such that \(\left \langle f, h \right \rangle=0\;\forall f \in M_1\) but \(\left \langle \varphi^{|\la|}_k,h \right \rangle \neq 0\). Since, \(\varphi^{|\la|}_k(-w)=\varphi^{|\la|}_k(w)\) and \(f(-w) \in M_1\) whenever \(f \in M_1\), we can assume that \(h(-w)=h(w)\). Thus, 
 % \begin{align*}
%(f\times_{-|\la|}\overline{h})(z)&=\int\limits_{\mathbb{C}^n}f(z-w)e^{\frac{-i}{2}|\la|\operatorname{Im}(z \cdot \overline{w})} \overline{h(w)}dw\\
%&=\int\limits_{\mathbb{C}^n}f(z+w)e^{\frac{i}{2}|\la|\operatorname{Im}(z \cdot \overline{w})} \overline{h(-w)}dw\\
%&=0.
%\end{align*}
%\textcolor{red}{contradiction???}
Now, we show that \Be \overline{\Phi^{\la}_{\al' ,\beta'}(w')\Phi^{\la}_{\al'' ,\beta''}(w'')} \in M_1, \Ee for \(|\al'|=j\) and \(|\al''|=k\). Let \(V\) be the orthogonal complement of \(M_1\) and \(g \in V\). Note that \(g(-w) \in V\) as well. It follows that \(\psi \times_{-\boldsymbol{\lambda}} \overline{g}=0\) which is same as \(g\times_{-\boldsymbol{\lambda}}\psi=0\). Taking  Weyl transform, we get 
$$
W_{-\boldsymbol{\lambda}}(g)W_{-\boldsymbol{\lambda}}(\psi)=0.
$$
Observe that 
\begin{align*}   \psi(w',w'')&=\sum_{|\gamma'|=j}\Phi^{\la}_{\gamma',\gamma'}(w')\sum_{|\gamma''|=k}\Phi^{\la}_{\gamma'',\gamma''}(w'')\\
&=\sum\limits_{\substack{(\gamma',\gamma'') \\ |\gamma'|=j,|\gamma''|=k}} \Phi^{\la}_{\gamma',\gamma'}(w')\Phi^{\la}_{\gamma'',\gamma''}(w'')\\
&=\sum\limits_{\substack{(\gamma',\gamma'') \\ |\gamma'|=j,|\gamma''|=k}} \Phi^{\la}_{(\gamma',\gamma''),(\gamma',\gamma'')}(w',w'')
\end{align*}
Using this in $$
\left \langle W_{-\boldsymbol{\lambda}}(g)W_{-\boldsymbol{\lambda}}(\psi)\Phi^{\la}_{(\al',\al'')},\Phi^{\la}_{(\beta',\beta'')} \right \rangle=0, \text{ for } |\al'|=j \text{ and } |\al''|=k,
$$
we obtain that $$
\int\limits_{\mathbb{C}^n}g(w)\Phi^{\la}_{(\al',\al''),(\beta',\beta'')}(w)dw=0.
$$
Note that the above computation involving the Weyl transform is similar to the case of the Heisenberg group, so we skip the computations here.
This implies that \Be 
\overline{\Phi^{\la}_{(\al',\al''),(\beta',\beta'')}(w)} = \overline{\Phi^{\la}_{\al',\beta'}(w')\Phi^{\la}_{(\al'',\beta'')}(w'')}\in M_1, 
\Ee for \(|\al'|=j\), \(|\al''|=k\). Thus implies that \(M\) is the entire space \(\mathcal{H}^{\boldsymbol{\lambda}}_{j,k}\).     
\end{proof}
%From now on, for identification of the Lie algebra of \(N\) with \(\mathbb{C}^n \times \mathbb{R}^3\), we consider the basis corresponding to \(\lambda\). 
%Let \(\varphi^{|\la|,2q}_{k'}(w')\) and \(\varphi^{|\la|,2l}_{k''}(w'')\) be the Laguerre function on \(\mathbb{C}^{2q}\)  and on \(\mathbb{C}^{2l}\) respectively.
Consider the function \Be
e^{\boldsymbol{\lambda}}_{j,k}(w,s)=e^{-i \boldsymbol{\lambda} \cdot s}\varphi^{2q-1}_{j,\la}(w')\varphi^{2l-1}_{k,\la}(w'')=e^{-i \boldsymbol{\lambda} \cdot s}\varphi^{q,l}_{j,k,\la}(w',w''),
\Ee
for $w=(w', w'')\in \mathbb C^{2q + 2l}, s\in \mathbb R^3$.
Then we have
\Be
\rho^{\boldsymbol{\lambda}}_{j,k}\left((\sigma',\sigma''),0,0\right)e^{\boldsymbol{\lambda}}_{j,k}(w,s)=e^{\boldsymbol{\lambda}}_{j,k}\left((\sigma',\sigma'')^{-1}w,s\right)=e^\la_{j,k}(w,s).
\Ee 
Therefore, the function \(e^{\boldsymbol{\lambda}}_{j,k}\) is an \(U\)-fixed vector. Since, \((U,N)\) forms a Gelfafnd pair, upto scalar multiple, \(e^{\boldsymbol{\lambda}}_{j,k}\) is the unique fixed vector.
Furthermore, for an integrable function \(f\) on \(N\), and function \(\psi(w,s)=e^{-i \boldsymbol{\lambda} \cdot s} \phi(w) \in \mathcal{H}^{\boldsymbol{\lambda}}_{j,k}\) we obtain that
\begin{align*}
\rho^{\boldsymbol{\lambda}}_{j,k}(f)\psi(w,s)&=\int\limits_{U \ltimes N}f(z,t)\rho^{\boldsymbol{\lambda}}_{j, k}(\sigma,z,t)\psi(w,s) \;d\sigma ds dt\\
&=\int\limits_{U \ltimes N}f(z,t)\psi(\sigma^{-1}w-\sigma^{-1}z,s-t+\frac{1}{2}[-\sigma^{-1}z,\sigma^{-1}w])\;d \sigma dz dt\\
&=\int\limits_{U \ltimes N}f(z,t)e^{-i\boldsymbol{\lambda} \cdot \left (-t+s+\frac{1}{2}[-z,w]\right )}\phi(-\sigma^{-1}z+\sigma^{-1}w)\;d\sigma dz dt\\
&=e^{-i\boldsymbol{\lambda} \cdot s}\int\limits_{U \times \mathbb{C}^{n}}f^{\boldsymbol{\lambda}}(z)e^{\frac{i}{2}\la \operatorname{Im}(z \cdot \overline{w})}\phi(\sigma^{-1}(w-z))\;d\sigma dz\\
&=e^{-i\boldsymbol{\lambda} \cdot s}(f^{\boldsymbol{\lambda}} \times_{\la}\Phi)(w),
\end{align*}
 where \Be 
 \Phi(w)=\displaystyle\int\limits_{Sp(q) \times Sp(l)}\phi(\sigma w)d\sigma.
 \Ee
 It can be proved that \(\rho^{\boldsymbol{\lambda}}_{j, k}(f)\) is a rank 1 operator. The action of the operator \(\rho^{\boldsymbol{\lambda}}_{j, k}(f)\) on the space \(\mathcal{H}^{\boldsymbol{\lambda}}_{j,k}\) can be summarised as
   $$\rho^{\boldsymbol{\lambda}}_{j,k}(f)(E^{\boldsymbol{\lambda}}_{\al,\beta})= \begin{cases}
    0\; &\textrm{if}\; \al=(\al',\al'') \neq (\beta',\beta'')=\beta \\
    (2\pi)^{-n/2}\frac{j!(2q-1)!}{(j+2q-1)!} \times \frac{k!(2l-1)!}{(k+2l-1)!}\left(f * e^{\boldsymbol{\lambda}}_{j,k}\right)\; &\textrm{if}\; \al=(\al',\al'')=(\beta',\beta'') =\beta.
\end{cases}  $$
The proof of the above is similar to the case of the Heisenberg motion group \cite{thangaveluheisenberg} and hence we omit it.
From this fact, we obtain that 
\Be
\rho^{\boldsymbol{\lambda}}_{j,k}(f)e^{\boldsymbol{\lambda}}_{j,k}=f \ast e^{\boldsymbol{\lambda}}_{j,k}.
\Ee
We have the following Fourier inversion formula for the Schwartz class function \(f\) on \(N\).  
\begin{prop}\label{inversion m=3}
    For every Schwartz class function f on \(N\), the following inversion formula holds
    $$
    f(z,t)=(2\pi)^{-n}\int\limits_{\mathbb{R}^3}{\sum\limits_{j,k}\left \langle \rho^{\boldsymbol{\lambda}}_{j,k}(f)e^{\boldsymbol{\lambda}}_{j,k}, \rho^{\boldsymbol{\lambda}}_{j,k}(1,z,t)e^{\boldsymbol{\lambda}}_{j,k} \right \rangle}\;d\mu(\boldsymbol{\lambda}),
    $$
    where \(n=2q+2l\).
\end{prop}
\begin{proof}
    Recall the following Fourier inversion formula for the H-type group \cite{narayananinjectivity},
    \begin{equation}
   f(z,t)=\int\limits_{\mathbb{R}^3}{\sum\limits_{r=0}^{\infty}f \ast e^{\boldsymbol{\lambda}}_r}\;d\mu(\boldsymbol{\lambda}).   
  \end{equation}
  Note that here, 
  $$e^{\boldsymbol{\lambda}}_r(w,s)=e^{-i \boldsymbol{\lambda} \cdot s}\varphi^{2n-1}_{r,\lambda}(w)=e^{-i \boldsymbol{\lambda} \cdot s}\sum\limits_{j+k=r}\varphi^{2q-1}_{j,\lambda}(w')\varphi^{2l-1}_{k,\lambda}(w'')=e^{-i \boldsymbol{\lambda} \cdot s}\sum\limits_{j+k=r}\varphi^{q,l}_{j,k,\lambda}(w',w'').$$
Thus, we can write the inversion formula as  
   \begin{equation}
   f(z,t)=\int\limits_{\mathbb{R}^3}\sum\limits_{r=0}^{\infty}\sum\limits_{j+k=r} f \ast e^{\boldsymbol{\lambda}}_{j,k}\;d\mu(\boldsymbol{\lambda})=\int\limits_{\mathbb{R}^3}\sum\limits_{j,k} f \ast e^{\boldsymbol{\lambda}}_{j,k}\;d\mu(\boldsymbol{\lambda}).   
  \end{equation}
Hence, it is enough to show that 
\Be
\left \langle \rho^{\boldsymbol{\lambda}}_{j,k}(f)e^{\boldsymbol{\lambda}}_{j,k}, \rho^{\boldsymbol{\lambda}}_{j,k}(1,z,t)e^{\boldsymbol{\lambda}}_{j,k} \right \rangle=(2\pi)^{n}f \ast e^{\boldsymbol{\lambda}}_{j,k}.
\Ee 
%The calculations are identical to those in the proof of Theorem \ref{inversion-formula}, so we will omit them.
For \(\al=(\al',\al'')\) with \(|\al'|=j\) and \(|\al''|=k\), we have \Be
\rho^{\boldsymbol{\lambda}}_{j,k}(1,z,t)E^{\boldsymbol{\lambda}}_{\al ,\al}((w,s)=E^{\boldsymbol{\lambda}}_{\al, \al}((-z,-t)(w,s)).
\Ee
\begin{comment}
Now, 
\begin{align*}
    E^{\boldsymbol{\lambda}}_{\al, \al}((-z,-t)(w,s))&=(2\pi)^{-n/2}\left \langle \pi_{-\boldsymbol{\lambda}}((-z,-t)(w,s)) \Phi^{\la}_\al, \Phi^{\la}_\al \right \rangle\\
    &=(2\pi)^{-n/2}\left \langle \pi_{-\boldsymbol{\lambda}}(w,s) \Phi^{\la}_{(\al',\al'')},\pi_{-\boldsymbol{\lambda}}(z,t) \Phi^{\la}_{(\al'\al'')} \right \rangle\\
    &=(2\pi)^{-n/2}\sum\limits_{(\beta',\beta'')}\left \langle \pi_{-\boldsymbol{\lambda}}(w,s)\Phi^{\la}_{(\al',\al'')}, \Phi^{\la}_{(\beta',\beta'')} \right \rangle\; \overline{\left \langle \pi_{-\boldsymbol{\lambda}}(z,t)\Phi^{\la}_{(\al',\al'')}, \Phi^{\la}_{(\beta',\beta'')} \right \rangle}\\
    &=(2\pi)^{n/2}\sum\limits_{\beta=(\beta',\beta'')}E^{\boldsymbol{\lambda}}_{\al, \beta}(w,s) \overline{E^{\boldsymbol{\lambda}}_{\al ,\beta}(z,t)}.
\end{align*}
\end{comment}
It is easy to verify that 
$$
E^{\boldsymbol{\lambda}}_{\al, \al}((-z,-t)(w,s))=(2\pi)^{n/2}\sum\limits_{\beta=(\beta',\beta'')}E^{\boldsymbol{\lambda}}_{\al, \beta}(w,s) \overline{E^{\boldsymbol{\lambda}}_{\al ,\beta}(z,t)}.
$$
Next, \(e^{\boldsymbol{\lambda}}_{j,k}=\displaystyle\sum\limits_{\substack{\al=(\al',\al'') \\ |\al'|=j,|\al''|=k}}(2\pi)^{n/2}E^{\boldsymbol{\lambda}}_{\al ,\al}\) , we get
\begin{align*}
    \left \langle \rho^{\boldsymbol{\lambda}}_{j,k}(f)e^{\boldsymbol{\lambda}}_{j,k}, \rho^\la_{j,k}(1,z,t)e^{\boldsymbol{\lambda}}_{j,k}\right \rangle&=\left \langle f \ast e^{\boldsymbol{\lambda}}_{j,k}, \displaystyle\sum\limits_{\substack{\al=(\al',\al'') \\ |\al'|=j,|\al''|=k}} \displaystyle\sum\limits_{\substack{\beta=(\beta',\beta'')}}(2\pi)^{n}E^{\boldsymbol{\lambda}}_{\al, \beta} \overline{E^{\boldsymbol{\lambda}}_{\al, \beta}} \;\right \rangle\\
    &=(2\pi)^{n}\displaystyle\sum\limits_{\substack{\al=(\al',\al'') \\ |\al'|=j,|\al''|=k}}\displaystyle\sum\limits_{\substack{\beta=(\beta',\beta'')}}\left \langle f*e^{\boldsymbol{\lambda}}_{j,k},E^{\boldsymbol{\lambda}}_{\al, \beta} \right \rangle E^{\boldsymbol{\lambda}}_{\al, \beta}(z,t) \\
    &=(2\pi)^{n} (f* e^{\boldsymbol{\lambda}}_{j,k})
\end{align*}
\end{proof}
Now we are in a position to state and prove Gutzmer's theorem for H type group $N$.
%When \(m=3\), the elements of N elements
 %can be written as \((z,w,t)\) with \(z\in \mathbbmss{H}^q \cong \mathbb{C}^{2q}, w \in \mathbbmss{H}^l \cong \mathbb{C}^{2l}, t \in \mathbb{R}^3 \), with \(2q+2l=n\).
%We identify \(N\) with \(\mathbb{R}^n\times \mathbb{R}^n \times \mathbb{R}^3\) and write its elements as \((x,u,t)\), where \(x=(x',x''),u=(u',u'')\) with \(x',u' \in \mathbb{R}^{2q}\) and \(x'',u'' \in \mathbb{R}^{2l}\). 
For a Schwartz class function \(f\) on \(N\), we say its Fourier transform is {\em compactly supported} if there exist two constants \(P\) and \(Q\) such that \Be f^{\boldsymbol{\lambda}}=0 \text{ for  } |\boldsymbol{\lambda}|>P
\Ee  and \Be \rho^{\boldsymbol{\lambda}}_{j,k}(f)=0 \text{ for } |\boldsymbol{\lambda}|(2(j+k)+n)>Q.\Ee
Therefore, for a Schwartz class function \(f\) with compact Fourier support, and for any \(g \in U \ltimes N\) we have 
$$
f(g.(x,u,\xi))=(2\pi)^{-n}\int\limits_{|\boldsymbol{\lambda}| \leq P}{\sum\limits_{|\boldsymbol{\lambda}|(2(j+k)+n) \leq Q} \left \langle \rho^{\boldsymbol{\lambda}}_{j,k}(f)e^{\boldsymbol{\lambda}}_{j,k},\rho^\la_{j,k}(g) \rho^{\boldsymbol{\lambda}}_{j,k}(1,x,u,\xi)e^{\boldsymbol{\lambda}}_{j,k} \right \rangle\; d\mu(\boldsymbol{\lambda})}.
$$
Furthermore, since each \(\rho^{\boldsymbol{\lambda}}_{j,k}(1,x,u,\xi)e^{\boldsymbol{\lambda}}_{j,k}\) extends to \(\mathbb{C}^{2n+3}\) as an entire function, so does the function \(f(g.(x,u,\xi))\). In this context, we state and prove an analogue of Gutzmer's formula for the H-type group $N$ with dimension of the center $m=3$. %The proof follows a similar framework to that presented in the preceding case, so we will omit the detailed justification for brevity.
\begin{theorem}\label{theorem 4.3}
    Let \(N\) denote the H-type group with center \(\mathbb{R}^3\) and \(f\) be a Schwartz class function on \(N\) whose Fourier transform is compactly supported. Then the function \(f\) extends to an entire function on \(\mathbb{C}^{2n+3}\) and we have the following estimate 
\begin{align*}
& \int_{U \ltimes N}|f(g.(z, w, \zeta))|^2 d g= \\
& (2 \pi)^{-n } \int\limits_{\mathbb{R}^3} e^{-2\boldsymbol{\lambda} \cdot \tau} e^{-\lambda (v \cdot x-u \cdot y)}\left(\sum_{j,k} C_{j,k} \left\|f^{\boldsymbol{\lambda}} \times_{\lambda} \varphi_{j,k,\lambda}^{q,l} \right\|_{L^2(\mathbb C^n)}^2  \varphi_{j,k,\lambda}^{q,l}(-2 i y, -2 i v)\right) d \mu(\boldsymbol{\lambda}),
\end{align*}    
where \(n=2q+2l\), \(|\boldsymbol{\lambda}|=\lambda ,\;z=x+iy, w=u+iv, \zeta=\omega+i\tau, y=(y',y''),v=(v',v'')\) and \(C_{j,k}=\frac{j!(2q-1)!}{(j+2q-1)!}\frac{k!(2l-1)!}{(k+2l-1)!}\).
\end{theorem}
\begin{proof}
Let \(G=U \ltimes N\) and \(\Lambda={\mathbb{R}^3}\setminus\{0\} \times \mathbb{N} \times \mathbb{N}\) in \cite[Proposition 4.3]{gutzmer}. For \(\delta=(\boldsymbol{\lambda},j,k) \in \Lambda\), we have the representation \(\rho^{\boldsymbol{\lambda}}_{j,k}\) of \(U \ltimes N\). Take \(\xi(\delta)=\rho^{\boldsymbol{\lambda}}_{j,k}(f)e^{\boldsymbol{\lambda}}_{j,k}\) and \(\eta(\delta)=\rho^{\boldsymbol{\lambda}}_{j,k}(1,z,w,\zeta)e^{\boldsymbol{\lambda}}_{j,k}\). % where \(z=x+iy,w=u+iv,\zeta=\omega+i\tau, x=(x',x''),y=(y',y'')\). 
To prove the required estimate, we will evaluate \(\left \| \rho^{\boldsymbol{\lambda}}_{j,k}(1,z,w,\zeta )e^{\boldsymbol{\lambda}}_{j,k} \right \|^2\) and use \cite[Proposition 4.3]{gutzmer}.

Note that \Be e^{\boldsymbol{\lambda}}_{j,k}(x_1,u_1,t_1)=e^{-i \boldsymbol{\lambda} \cdot t_1 } \varphi^{q,l}_{j,k}(x_1,u_1)=e^{-i \boldsymbol{\lambda} \cdot t_1 }\varphi^{2q-1}_{j,\lambda}(x'_1,u'_1)\varphi^{2l-1}_{k,\lambda}(x''_1,u''_1),\Ee where \(x_1=(x'_1,x''_1)\in \mathbb{R}^{2q}\times \mathbb{R}^{2l}\) and \(u_1=(u'_1,u''_1)\in \mathbb{R}^{2q}\times \mathbb{R}^{2l}\), can be extended holomorphically on \(\mathbb{C}^n \times \mathbb{C}^n \times \mathbb{C}^m\). %Moreover, this function can be extended to \(\mathbb{C}^{2n+m}\).
We have the action 
\begin{align*}
   &\rho^{\boldsymbol{\lambda}}_{j,k}(1,x,u,t)e^{\boldsymbol{\lambda}}_{j,k}(x_1,u_1,t_1)\\&=e^{\boldsymbol{\lambda}}_{j,k}((1,x,u,t)^{-1}(x_1,u_1,t_1))\\
  % &=e^\la_{k',k''}((x,u,t)^{-1}(x_1,u_1,t_1) )\\
   %&=e^\la_{k',k''}((x_1-x,u_1-u),t_1-t+\frac{1}{2}[-z',w'])\;\;\text{where } z'=(x,u)\;,w'=(x_1,u_1)\\
  % &=e^{-i  \la \cdot ([-z',w'])}e^{-i \la \cdot (t_1-t)}\varphi^{|\la|,2q}_{k'}(x'_1-x',u'_1-u')\varphi^{|\la|,2l}_{k''}(x''_1-x'',u''_1-u'')\\
   %&=e^{i |\la|\operatorname{Im}(z' \cdot \overline{w'})}e^{-i \la \cdot (t_1-t)}\varphi^{|\la|}_k(x_1-x,u_1-u)\\
   &=e^{i\frac{\la}{2}(-x\cdot u_1+u\cdot x_1)}e^{-i \boldsymbol{\lambda} \cdot (t_1-t)}\varphi^{2q-1}_{j,\lambda}(x'_1-x',u'_1-u')\varphi^{2l-1}_{k,\lambda}(x''_1-x'',u''_1-u'').
\end{align*}
%The function \((x,u,t) \mapsto \rho^\la_{k',k''}(1,x,u,t)e^\la_{k',k''}\) can be extended holomorphically to \(\mathbb{C}^{2n+m}\) as 
This action has a holomorphic extension as
\begin{align*}
\rho^{\boldsymbol{\lambda}}_{j,k}(1,z,w,\zeta)e^{\boldsymbol{\lambda}}_{j,k}(x_1,u_1,t_1)=e^{-i \boldsymbol{\lambda} \cdot (t_1-\zeta )}e^{\frac{\la}{2}(w\cdot x_1-z\cdot u_1)}\varphi^{2q-1}_{j,\lambda}(x'_1-z',u'_1-w')\\
\varphi^{2l-1}_{k,\lambda}(x''_1-z'',u''_1-w''),    
\end{align*}
where \(z=(z',z'')\in \mathbb{C}^{2q}\times \mathbb{C}^{2l},w=(w',w'')\in \mathbb{C}^{2q} \times \mathbb{C}^{2l}\) and \(\zeta \in \mathbb{C}^3\).

Writing the above in terms of real and imaginary parts with \(z=x+iy\), \(w=u+iv\), and \(\zeta=\omega+i\tau\), we obtain
\begin{align*}
  \phantom{{\rho}^\la_k(1,z,w,\zeta)e^{\lambda}_k(x_1,u_1,t_1)}
  &\begin{aligned}
    \mathllap{{\rho}^{\boldsymbol{\lambda}}_{j,k}(1,z,w,\zeta)e^{\boldsymbol{\lambda}}_{j,k}(x_1,u_1,t_1)} &=e^{-\boldsymbol{\lambda} \cdot \tau}e^{-i \boldsymbol{\lambda} \cdot ( t_1-\omega) }e^{i\frac{\la}{2}(u\cdot x_1-x\cdot u_1)}e^{-\frac{\la}{2}(v\cdot x_1-y\cdot u_1)}\\
      &\qquad \qquad  \varphi^{2q-1}_{j,\lambda}(x'_1-z',u'_1-w')\varphi^{2l-1}_{k,\lambda}(x''_1-z'',u''_1-w'').
      %e^{-i \la \cdot \left ( t_1-(\omega+i\tau)\right )}e^{|\la|((u+iv)\cdot x_1-(x+iy)\cdot u_1)} \\
     % &\qquad \qquad \varphi^{|\la|,2q}_{k'}(x'_1-z',u'_1-w')\varphi^{|\la|,2l}_{k''}(x''_1-z'',u''_1-w'')
  \end{aligned}\\
  %&\begin{aligned}
   % \mathllap{\;} &= e^{\la \cdot \tau}e^{-i \la \cdot ( t_1-\omega) }e^{i|\la|(u\cdot x_1-x\cdot u_1)}e^{-\frac{|\la|}{2}(v\cdot x_1-y\cdot u_1)}\\
     % &\qquad \qquad  \varphi^{|\la|,2q}_{k'}(x'_1-z',u'_1-w')\varphi^{|\la|,2l}_{k''}(x''_1-z'',u''_1-w'').
  %\end{aligned}
\end{align*}
Therefore we have,
\begin{align}
    &{ \left \| \rho^{\boldsymbol{\lambda}}_{j,k}(1,z,w,\zeta)e^{\boldsymbol{\lambda}}_{{j,k}}  \right \|}^2 \nonumber \\ 
 &={(\la)}^n\int\limits_{\mathbb{R}^n \times \mathbb{R}^n}{e^{-2 \boldsymbol{\lambda} \cdot \tau}e^{-\la(v\cdot x_1-y\cdot u_1)}{|\varphi^{2q-1}_{j,\lambda}(x'_1-z',u'_1-w')|}^2{|\varphi^{2l-1}_{k,\lambda}(x''_1-z'',u''_1-w'')|}^2\;dx_1du_1}  \nonumber \\
  &=e^{-2 \boldsymbol{\lambda} \cdot \tau} I_1 \times I_2, \label{I_ times I_2} 
\end{align}
    where 
    \begin{equation*}
    I_1=\int\limits_{\mathbb{R}^{2q} \times \mathbb{R}^{2q}}{e^{-(\la)^{1/2}(v' \cdot x'_1-y' \cdot u'_1)}\,{|\varphi^{2q-1}_{j}(x'_1-{(\la)}^{1/2}z',u'_1-{(\la)}^{1/2}w')|}^2\;dx'_1du'_1},
    \end{equation*} and 
    \begin{equation*}
        I_2=\int\limits_{\mathbb{R}^{2l} \times \mathbb{R}^{2l}}{e^{-(\la)^{1/2}(v'' \cdot x''_1-y'' \cdot u''_1)}\,{|\varphi^{2l-1}_{k}(x''_1-{(\la)}^{1/2}z'',u''_1-{(\la)}^{1/2}w'')|}^2\;dx''_1du''_1}.
    \end{equation*}
    Now proceeding as in proof of \cite[ Theorem 4.2]{gutzmer}, we obtain
    \begin{equation*}
        I_1=e^{-(\la)(v' \cdot x'-y' \cdot u')}\varphi^{2q-1}_{j}(-2i{(\la)}^{1/2}y',-2i{(\la)}^{1/2}v'),
    \end{equation*} and 
    \begin{equation*}
         I_2=e^{-(\la)(v'' \cdot x''-y'' \cdot u'')}\varphi^{2l-1}_{k}(-2i{(\la)}^{1/2}y'',-2i{(\la)}^{1/2}v'').
    \end{equation*}
    Substituting these in Eq. \eqref{I_ times I_2}, we deduce that 
    \begin{align*}
    {\left \| \rho^{\boldsymbol{\lambda}}_{j,k}(1,z,w,\zeta)e^{\boldsymbol{\lambda}}_{j,k}  \right \|}^2&= e^{-2 \boldsymbol{\lambda} \cdot \tau}e^{-\la(v \cdot x-y \cdot u)}\varphi_{j,\lambda}^{2q-1}(-2 i y', -2 i v')\varphi_{k,\lambda}^{2l-1}(-2 i y'', -2 i v'')\\
    &=e^{-2 \boldsymbol{\lambda} \cdot \tau}e^{-\la(v \cdot x-y \cdot u)}\varphi^{q,l}_{j,k,\lambda}(-2iy,-2iv).
    \end{align*}
    This completes the proof of the theorem.
\end{proof}
\subsection{The case when $m=2$:}
In this case, \(N \cong \mathbb{C}^{2q}\) and \(U=Sp(q)\) acts transitively on the spheres centered at origin. Aside from a few technical details, this situation is quite similar to that of the Heisenberg group. To ensure that this paper is self-contained, we will briefly outline the proof, highlighting differences from the Heisenberg group case.

For \(k \in \mathbb{N},\; 0 \neq \boldsymbol{\lambda}=\lambda \omega \in \mathbb{R}^2\), let \(\mathcal{H}^{\boldsymbol{\lambda}}_k\) be the Hilbert space with orthonormal basis 
\begin{align*}
 E^{\boldsymbol{\lambda}}_{\al,\beta}(x,u,t)=%&(2\pi)^{-n/2} \left \langle \pi_\la(x,u,t)\Phi^{|\la|}_\al, \Phi^{|\la|}_\beta \right \rangle\\
 %&=\;e^{i\la \cdot t}(2\pi)^{-n/2} \left \langle \pi_\la(x,u,0)\Phi^{|\la|}_\al, \Phi^{|\la|}_\beta  \right \rangle\\
 %&=e^{i\la \cdot t}(2\pi)^{-n/2} \left \langle \rho_{|\la|}(x,u,0)\Phi^{|\la|}_\al, \Phi^{|\la|}_\beta  \right \rangle\\
 &=e^{i\boldsymbol{\lambda} \cdot t}\Phi^{\la}_{\al ,\beta}(x,u),
\end{align*}%$$\varphi_k^{|\la|} (x, u)=\varphi_k\left(|\la|^{1 / 2}(x, u)\right), \text{ where }\varphi_k \text{ is the Laguerre function}$$
%and $$ e_k^{\boldsymbol{\lambda}}(x,u,t)=e^{-i\boldsymbol{\lambda} \cdot t} \varphi_{k,\la}^{2q-1}(x,u),$$
where \(|\boldsymbol{\lambda}|=\lambda,\; \al,\beta \in \mathbb{N}^{2q}\) with \(|\al|=k\).
We equip the space \(\mathcal{H}^{\boldsymbol{\lambda}}_k\) with the following inner product 
$$
\left \langle F,G\right \rangle= {( \la )}^{2q} \int\limits_{\mathbb{C}^{2q}}{F(z,0)\;\overline{G(z,0)} \;dz}.
$$
%For \(0 \neq \lambda \in \mathbb{R}^2 \) and \(k \in \mathbb{N}\), let \(\mathcal{H}^\la_k\) be the Hilbert space with orthonormal basis \(E^{-\la}_{\al,\beta}\) where \(\alpha, \beta \in \mathbb{N}^n\) and \(\lvert \alpha \rvert =k\). 
Moreover, elements of the space \(\mathcal{H}^{\boldsymbol{\lambda}}_k\) can be characterised by the following eigenvalue equations 
$$
\mathcal{L}(f)=(2k+n)\lvert \boldsymbol{\lambda} \rvert f=(2k+n)\la f, \;\;\;i\frac{\partial}{\partial t_j}=\la_j f, \;\; 1 \leq j \leq 2,
$$
where \(n=2q\).
%We equip the space \(\mathcal{H}^\la_k\) with the following inner product 
%$$
%\left \langle F,G\right \rangle= {\lvert \la \rvert}^n \int\limits_{\mathbb{C}^n}{F(z,0)\;\overline{G(z,0)} \;dz}.
%$$
Define a representation \(\rho^{\boldsymbol{\lambda}}_k \) of \(U \ltimes N\) on \(\mathcal{H}^{\boldsymbol{\lambda}}_k\) as
$$
\rho^{\boldsymbol{\lambda}}_k(\sigma,z,t)\varphi(w,s)=\varphi \left( (\sigma,z,t)^{-1}(w,s)\right), \;\;\; \text{where}\; \varphi \in \mathcal{H}^{\boldsymbol{\lambda}}_k.
$$
For the representation defined above, we obtain the following result. 
\begin{theorem}\label{irreducibility proof}
The representation \(\rho^{\boldsymbol{\lambda}}_k \) as defined above is an irreducible unitary representation of the group \(U \ltimes N\).
\end{theorem}
\begin{proof}
The proof is similar to the case when \(m=3\).
\end{proof}
For the function \(e^{\boldsymbol{\lambda}}_k(w,s)=e^{-i \boldsymbol{\lambda} \cdot s}\varphi^{2q-1}_{k,\lambda}(w)\), we have 
$$
\rho^{\boldsymbol{\lambda}}_k(\sigma,0,0)e^{\boldsymbol{\lambda}}_k(w,s)=e^{\boldsymbol{\lambda}}_k(\sigma^{-1}w,s)=e^{\boldsymbol{\lambda}}_k(w,s).
$$ 
Hence, it is a \(U\)-fixed vector. Since \((U,N)\) forms a Gelfand pair, upto scalar multiple, \(e^{\boldsymbol{\lambda}}_k\) is the unique fixed vector.

Given an integrable function \(f\) on \(N\), it can be considered as \(U\) invariant function on \(U \ltimes N\). Thus, the Fourier transform of function f can be defined conventionally as the operator 
$$
\rho^{\boldsymbol{\lambda}}_k(f)=\int\limits_{U \ltimes N}{f(g)\rho^{\boldsymbol{\lambda}}_k(g)dg}.
$$
\begin{prop}\label{inversion-formula}
    For every Schwartz class function f on \(N\), the following inversion formula holds
    $$
    f(z,t)=(2\pi)^{-n}\int\limits_{\mathbb{R}^2}{\sum\limits_{k=0}^{\infty}\left \langle \rho^{\boldsymbol{\lambda}}_k(f)e^{\boldsymbol{\lambda}}_k, \rho^{\boldsymbol{\lambda}}_k(1,z,t)e^{\boldsymbol{\lambda}}_k \right \rangle}\;d\mu(\boldsymbol{\lambda}),
    $$
    where \(n=2q\).
\end{prop}

%This concludes the representation theory and Fourier inversion formula needed to establish Gutzmer's formula for the case when \(m=2\). 
For a Schwartz class function \(f\) on \(N\), we say its Fourier transform is {\em compactly supported} if there exist two constants \(P\) and \(Q\) such that \Be f^{\boldsymbol{\lambda}}=0 \text{ for  } |\boldsymbol{\lambda}|>P
\Ee  and \Be \rho^{\boldsymbol{\lambda}}_{k}(f)=0 \text{ for } |\boldsymbol{\lambda}|(2k+n)> Q.\Ee Therefore, for a Schwartz class function with compact Fourier support, we have 
$$
f(g.(x,u,\xi))=(2\pi)^{-n}\int\limits_{|\boldsymbol{\lambda}| \leq P}{\sum\limits_{|\boldsymbol{\lambda}|(2k+n) \leq Q} \left \langle \rho^{\boldsymbol{\lambda}}_k(f)e^{\boldsymbol{\lambda}}_k,\rho^{\boldsymbol{\lambda}}_k(g) \rho^{\boldsymbol{\lambda}}_k(1,x,u,\xi)e^{\boldsymbol{\lambda}}_k \right \rangle\; d\mu(\boldsymbol{\lambda})}.
$$
Furthermore, since each \(\rho^{\boldsymbol{\lambda}}_k(1,x,u,\xi)e^{\boldsymbol{\lambda}}_k\) extends to \(\mathbb{C}^{2n+2}\) as an entire function, so does the function \(f(g.(x,u,\xi))\).
We obtain the following analogue of Gutzmer's formula for the action of H-type motion group on \(\mathbb{C}^{2n+2}\), whose proof is similar to the case when \(m=3\).
\begin{theorem}\label{theorem 4.2}
    Let \(N\) denote H-type group with center \(\mathbb{R}^2\) and \(f\) be a Schwartz class function on \(N\) whose Fourier transform is compactly supported. Then the function \(f\) extends to an entire function on \(\mathbb{C}^{2n+2}\) and we have the following estimate
   $$
\begin{aligned}
& \int_{U \ltimes N}|f(g.(z, w, \zeta))|^2 d g= \\
& (2 \pi)^{-n } \int\limits_{\mathbb{R}^2} e^{-2 \boldsymbol{\lambda} \cdot \tau} e^{- \lambda (v \cdot x-u \cdot y)}\left(\sum_{k=0}^{\infty}\left\|f^{\boldsymbol{\lambda}} \times_{ \lambda } \varphi_{k,\lambda}^{n-1}\right\|_2^2 \frac{k!(n-1)!}{(k+n-1)!} \varphi_{k,\lambda}^{n-1}(-2 i y, -2 i v)\right) d \mu(\boldsymbol{\lambda}),
\end{aligned}
$$    
where \(n=2q\), \(|\boldsymbol{\lambda}|=\lambda\), \(z=x+iy, w=u+iv, \zeta=\omega+i\tau\), $\left\|f^{\boldsymbol{\lambda}} \times_{\lambda} \varphi_{k,\lambda}^{2q-1}\right\|_2$ is the $L^2\left(\mathbb{C}^n\right)$ norm of $f^{\boldsymbol{\lambda}} \times_{\lambda} \varphi_{k,\lambda}^{2q-1}$.
\end{theorem}

\section{Gutzmer's formula on reduced H-type group}\label{sec5}
In this section, we define reduced H-type group and prove an analogue of Gutzmer's formula for the H-type group $N$. We need this Gutzmer's formula for the proof of Beurling's theorem. The theory of Fourier transform and Plancherel theorem for the case of reduced Heisenberg group has been covered in \cite{dasgupta-santosh-localization-operator}. 

Let \(N_1=\{(0,0,2 \pi t): t \in \mathbb{Z}^m\}\). It is easy to check that \(N_1\) is a normal subgroup of \(N\). Consider the quotient group \(N/N_1\) with the group operation 
\begin{equation}\label{group law on reduced}
    (z,t)(w,s)=\left(z+w,t+s+\frac{1}{2} [z,w] \right),
\end{equation}
where the coordinates of \(t+s+\frac{1}{2}[z,w]\) are modulo \(2\pi\).

We define the reduced H-type group \(N_{red}\) as the quotient group \(N/N_1\) with the above group law.\\
One can check that the Schr\"odinger representation \(\pi_{\boldsymbol{\lambda}}\) is a representation of \(N_{red}\) iff \(\boldsymbol{\lambda} \in \mathbb{Z}^m\).
\begin{defi}
    Let \(f \in L^1(N_{red})\) and \(\mathbf{p} \in \mathbb{Z}^m\setminus \{0\}\). The Fourier transform of \(f\) at \(\mathbf{p}\) is defined as
    \Be
    \Hat{f}(\mathbf{p})=\displaystyle\int\limits_{N_{red}}f(z,t)\pi_{\mathbf{p}}(z,t)\frac{dzdt}{(2\pi)^m}.
    \Ee
\end{defi}
For \(\mathbf{p} \in \mathbb{Z}^m\), let \(f^{\mathbf{p}}(z)=\displaystyle\int\limits_{\mathbb{T}^m}f(z,t) e^{i \mathbf{p}\cdot t}\frac{dt}{(2\pi)^m}\), be the inverse Fourier transform on \(\mathbb{T}^m\).
Then we have  \Be 
\Hat{f}(\mathbf{p})= \int_{\mathbb{C}^n}f^{\mathbf{p}}(z)\pi_{\mathbf{p}}(z) dz.
\Ee
We deduce the following Plancherel theorem for the reduced H-type group, whose proof is the same as the Heisenberg group case, hence, we are omitting the proof.
\begin{theorem}
    Let \(f \in L^2(N_{red})\) and \(\mathbf{p} \in \mathbb{Z}^m\setminus \{0\}\). Then \(\Hat{f}(\mathbf{p}): L^2(\mathbb{R}^n \to L^2(\mathbb{R}^n)\) is a Hilbert Schmidt operator and we have the following identity
    $$
    \sum_{\mathbf{p} \in {\mathbb{Z}^m}\setminus \{0\}} {\lVert \Hat{f}(\mathbf{p})\rVert}^2_{HS} = (2\pi)^{n+m}|\mathbf{p}|^{-n}{\lVert f-f^0 \rVert}_{2}^2.
    $$
  % where \(f^0(z)=\displaystyle\int\limits_{\mathbb{T}^m} f(z,t)\frac{dt}{(2\pi)^m}\).

\end{theorem}
We say a function \(f\) on \(N_{red}\) has {\em mean value \(0\)}  if \Be 
\int_{\mathbb{T}^m}f(z,t)dt=0.\Ee
Let \(U \ltimes N_{red}\) be the reduced H-type motion group. Then, following the same method as in the case of the H-type group, one can deduce Gutzmer's formula for the reduced H-type motion group.

\subsection{The case when $m=3$:}
For the case when the center of \(N\) is \(\mathbb{R}^3\), we have the following analogue of Gutzmer's formula for its corresponding reduced H-type group. The proof is on the same lines as the proof of Gutzmer's formula for the case of H-type group with \(m=3\).% The proof is on the same lines as Theorem \ref{Gutzmers formula for reduced} with minor adjustments. Thus, we are omitting its proof. 
\begin{theorem}
    Let \(N\) be the H-type group with center \(\mathbb{R}^3\) and \(N_{red}\) be its corresponding reduced H-type group. Let \(f\) be a function on \(N_{red}\) having mean value 0, satisfying the conditions stated in Theorem \ref{theorem 4.3}. Then we have 
    \small{
    \begin{align*}
   &\int_{U \ltimes N_{red}}|f(g.(z, w, \zeta))|^2 d g= \\
&(2 \pi)^{-2n-3 } \sum_{\mathbf{p} \neq 0} e^{-2 \mathbf{p} \cdot \tau} e^{- p (v \cdot x-u \cdot y)}\left(\sum_{j,k}C_{j,k}\left\|f^{\mathbf{p}} \times_{ p } \varphi_{j,k,p}^{q,l}\right\|_2^2  \varphi_{j,k,p}^{q,l}(-2 i y, -2 iv)\right)|\mathbf{p}|^n,
   \end{align*}
   }
   where \(n=2q+2l\), \(|\mathbf{p}|=p\), \(C_{j,k}=\frac{j!(2q-1)!}{(j+2q-1)!}\frac{k!(2l-1)!}{(k+2l-1)!}\), \(z=x+iy\), \(w=u+iv\), \(y=(y',y''),v=(v',v'')\) and \(\zeta=\omega+i\tau\).
\end{theorem}
From this, we have
\begin{align}
  \int_{U \ltimes N_{red}}|f(g.(iy, iv,t))|^2 d g&= 
(2 \pi)^{-2n-3 } \sum_{\mathbf{p} \neq 0} \sum_{j,k}C_{j,k}\left\|f^{\mathbf{p}} \times_{p} \varphi_{j,k,p}^{q,l}\right\|_2^2 \varphi_{j,k,p}^{q,l}(-2 i y, -2 iv)|\mathbf{p}|^n. \label{red gutzmer m=3 eq1}
\end{align}
Thus, in particular, for function of the form \Be 
f(x,y,t)=h\left( \frac{x}{\sqrt{\abs{\boldsymbol{\la}}}}, \frac{y}{\sqrt{\abs{\boldsymbol{\la}}}}\right)e^{-i (\mathbf{1}\cdot t )}=h\left( \frac{x}{\sqrt{\lambda}}, \frac{y}{\sqrt{\lambda}}\right)e^{-i (\mathbf{1}\cdot t )},\Ee  
where \(\boldsymbol{\la}=\lambda \omega \in \mathbb{R}^3 \setminus \{0\}\), with \(|\boldsymbol{\lambda}|=\lambda\) and \(\mathbf{1}=(1,1,1)\) we have \Be f^{\mathbf{p}}=0, \Ee except for \(\mathbf{p}=\mathbf{1}\). Thus, for such functions, Eq. \eqref{red gutzmer m=3 eq1} takes the following form
\begin{align}
\int_{N^3_{red}}|f(g.(iy, iv,t))|^2 d g&= 
(2 \pi)^{-2n-3 }{3}^{n/2} \sum_{j,k}C_{j,k}\left\|f^{\mathbf{1}} \times_{\sqrt{3}} \varphi_{j,k,\sqrt{3}}^{q,l}\right\|_2^2  \varphi_{j,k,\sqrt{3}}^{q,l}(-2 i y, -2 iv). \label{red gutzmer m=3 eq2}
\end{align}
Furthermore, if \(g=(\sigma,x,u,s)\), then 
\begin{align*}
    {\left| f(g.(iy,iv,t) \right |}^2&={\left | h \left ( \frac{1}{\sqrt{\la}} \left((x,u)+i\sigma(y,v) \right) \right ) e^{-i \left ( \mathbf{1}\cdot (s+t+\frac{i}{2}[((x,u),\sigma(y,v)]) \right )}\right |}^2\\
    &={\left | h \left ( \frac{1}{\sqrt{\la}} \left((x,u)+i\sigma(y,v) \right) \right ) \right |}^2 e^{-\left ( \mathbf{1}\cdot [(x,u),\sigma(y,v)] \right )}.
\end{align*}

Now, 
\begin{align*}
    &\left( f^{\mathbf{1}} \times_{\sqrt{3}} \varphi_{j,k,\sqrt{3}}^{q,l}\right )(x,u)\\
   &= \int_{\mathbb{R}^{2n}}  h \left ( \frac{1}{\sqrt{\lambda}} \left(x-y,u-v \right) \right ) \varphi_{j,k,\sqrt{3}}^{q,l}(y,v) e^{i\sqrt{3}\operatorname{Im}\left ((x,u) \cdot \overline{(y,v)}\right )} dydv\\
   &=(\la)^n \int_{\mathbb{R}^{2n}} h \left ( \frac{x}{\sqrt{\lambda}}-y,\frac{u}{\sqrt{\lambda}} -v\right )\varphi_{j,k,\sqrt{3}\lambda}^{q,l}(y,v)e^{i\sqrt{3}\la \operatorname{Im}\left ( \left ( \frac{x}{\sqrt{\lambda}},\frac{u}{\sqrt{\lambda}}\right ) \cdot \overline{\left( y,v\right )}\right )}dydv\\
   &=(\la)^n\left ( h \times_{\sqrt{3}\lambda}\varphi^{q,l}_{j,k,\sqrt{3}\lambda}\right ) \left( \frac{x}{\sqrt{\lambda}},\frac{u}{\sqrt{\lambda}}\right).
   \end{align*}
 Therefore,
\begin{align*}
    {\lVert f^1\times_{\sqrt{3}}\varphi_{j,k,\sqrt{3}}^{q,l}\rVert}^2_2&= \int_{\mathbb{R}^{2n}}(\la)^{2n}{\left | \left ( h \times_{\sqrt{3}\lambda}\varphi^{q,l}_{j,k,\sqrt{3}\lambda}\right ) \left( \frac{x}{\sqrt{\lambda}},\frac{u}{\sqrt{\lambda}}\right) \right |}^2dxdu\\
    &=(\la)^{3n}{\lVert h\times_{\sqrt{3}\lambda}\varphi^{q,l}_{j,k,\sqrt{3}\lambda}\rVert}^2_2.
\end{align*}
Thus, using Eq. \eqref{red gutzmer m=3 eq2} we obtain  
\begin{align*}
    &\int_{U}\int_{\mathbb{R}^{2n}}{\left | h\left ( (x,u)+i\sigma(y,v)\right ) \right |}^2e^{-\sqrt{3}\lambda \operatorname{Im}\left( (x,u) \cdot \overline{\sigma(y,v)}\right)} d\sigma dx du \nonumber\\
    &=(\la)^{2n}3^{n/2}(2\pi)^{-2n-3} \sum_{j,k} C_{j,k}{\left \lVert h\times_{\sqrt{3}\lambda}\varphi^{q,l}_{j,k,\sqrt{3}\lambda}\right \rVert}^2_2  \varphi^{q,l}_{j,k,\sqrt{3}\lambda}(-2iy,-2iv).
    \end{align*}
    That is, for each \(\boldsymbol{\la} \in \mathbb{R}^3\setminus \{0\}\) with \(|\boldsymbol{\lambda|}=\lambda\), we have the following 
    \begin{align}\label{red gutzmer eq3}
        &\int_{U}\int_{\mathbb{R}^{2n}}{\left | h\left ( (x,u)+i\sigma(y,v)\right ) \right |}^2e^{-\la \operatorname{Im}\left( (x,u) \cdot \overline{\sigma(y,v)}\right)} d\sigma dx du \nonumber\\
    &=(\la)^{2n}3^{-n/2}(2\pi)^{-2n-3} \sum_{j,k} C_{j,k}{\left \lVert h\times_{\lambda}\varphi^{q,l}_{j,k,\lambda}\right \rVert}^2_2  \varphi^{q,l}_{j,k,\lambda}(-2iy,-2iv) \nonumber \\
&=(\la)^{2n}3^{-n/2}(2\pi)^{-2n-3} \sum_{j,k} C_{j,k}{\left \lVert h\times_{\lambda}\varphi^{2q-1}_{j,\lambda}\varphi^{2l-1}_{k,\lambda}\right \rVert}^2_2 \varphi^{2q-1}_{j,\lambda}(-2iy',-2iv')\varphi^{2l-1}_{k,\lambda}(-2iy'',-2iv'').
    \end{align}
The action of \(U\) on \(\mathbb{R}^{2n}\) can be extended to \(\mathbb{C}^{2n}\). The next proposition plays a crucial role in proving Beurling's theorem. 
\begin{prop}\label{beurling estimate m=3}
    Let \(N\) be H-type group with center \(\mathbb{R}^3\) and \(f\) be a Schwartz class function on \(N\), Then we have the following estimate
    \small{$$ 
    \begin{aligned}
        &\displaystyle\int_{U}{\lVert \pi_{\boldsymbol{\lambda}}(\sigma(z,w))\Hat{f}(\boldsymbol{\lambda})\rVert}^2_{HS}d\sigma\\
        &=e^{-\frac{\la}{2}(v\cdot x-u\cdot y)}(2\pi)^{-n-3} \sum_{j,k}{\left( \frac{\la}{\sqrt{3}} \right )}^{n}C_{j,k}{\left \lVert f^{\boldsymbol{\lambda}}\times_{\lambda}\varphi^{2q-1}_{j,\lambda}\varphi^{2l-1}_{k,\lambda} \right \rVert}^2_2\varphi^{2q-1}_{j,\lambda}(2iy',2iv')\varphi^{2l-1}_{k,\lambda}(2iy'',2iv''),
    \end{aligned}
    $$}
    where \(|\boldsymbol{\lambda}|=\lambda\), \(C_{j,k}=\frac{j!(2q-1)!}{(j+2q-1)!}\frac{k!(2l-1)!}{(k+2l-1)!}\), \(z=x+iy\) and \(w=u+iv\). 
\end{prop}
\begin{proof}
Observe that, \(\displaystyle\int_U {\lVert \pi_{\boldsymbol{\lambda}}(\sigma(z,w))\Hat{f}(\boldsymbol{\lambda})\rVert}^2_{HS}d\sigma= \sum_{l',l''}\displaystyle\int_U{\lVert \pi_{\boldsymbol{\lambda}}(\sigma(z,w))\Hat{f}(\boldsymbol{\lambda})P_{l',l''}(\la)\rVert}^2_{HS}d\sigma\)\\
where \(P_{l',l''}(\la)\) is the projection on the space \(\left \{ \Phi_{\al'}^{\la }(z')\Phi_{\al''}^{\la}(z''): |\al'|=l', |\al''|=l'' \right \}.\)

Let \(\psi_{l',l''}(x,u)=f^{\boldsymbol{\lambda}}\times_{\la}\varphi^{2q-1}_{l',\lambda}\varphi^{2l-1}_{l'',\lambda}\) so that \(\Hat{f}(\boldsymbol{\lambda})P_{l',l''}(\la)=(2\pi)^{-n}|\boldsymbol{\lambda}|^n\pi_{\boldsymbol{\lambda}}(\psi)\). Note that the function \(\psi_{l',l''}\) admits a holomorphic extension.
Next, for \(\boldsymbol{\nu} \in \mathbb{R}^3\setminus \{0\}\) consider the twisted translation of \(\psi\) defined as 
\begin{align*}
   \tau^{\boldsymbol{\nu}}_{(x,u)}\psi_{l',l''}(x',u')&=\psi_{l',l''}(x'-x,u'-u)e^{i\frac{\boldsymbol{\nu}}{2} \cdot [(x,u),(x',u')]}\\
   &=\psi_{l',l''}(x'-x,u'-u)e^{i\frac{|\boldsymbol{\nu}|}{2}(-xu'+ux')}.
\end{align*}
It is a routine calculation to check that \(\pi_{\boldsymbol{\nu}}(x,u)\pi_{\boldsymbol{\nu}}\psi_{l',l''}=\pi_{\boldsymbol{\nu}}\left( \tau^{\boldsymbol{\nu}}_{(x,u)}\psi_{l',l''} \right )\).\\
Using Eq. \eqref{red gutzmer eq3}, we have,
\begin{align*}
    &\sum_{l',l''}\displaystyle\int_U{\lVert \pi_{\boldsymbol{\lambda}}(\sigma(z,w))\Hat{f}(\boldsymbol{\lambda})P_{l',l''}(\boldsymbol{\lambda})\rVert}^2_{HS}d\sigma\\
    &=\sum_{l',l''}\displaystyle\int_U(2\pi)^{-n}(\lambda)^{n}\int_{\mathbb{R}^{2n}}{\left | \tau^{\boldsymbol{\lambda}}_{\sigma(z,w)}\psi_{l',l''}(x',u')\right |}^2dx'du'\\
    &=\sum_{l',l''}(2\pi)^{-n}{(\lambda)}^ne^{-\frac{\lambda}{2}(v\cdot x-u\cdot y)}\int_U\int_{\mathbb{R}^{2n}}{\left | \psi_{l',l''} \left( (x',u')-i\sigma(y,v) \right)\right |}^2e^{-\frac{\lambda}{2}\operatorname{Im}\left( -(x',u') \cdot \overline{\sigma(y,v)}\right)}dx'du'd\sigma\\
    &=\sum_{l',l''}(2\pi)^{-n}(\lambda)^{3n}e^{-\frac{\lambda}{2}(v\cdot x-u\cdot y)}3^{-n/2}(2\pi)^{-2n-3} \sum_{j,k}{\left \lVert \psi_{l',l''} \times_{\lambda
    }\varphi^{2q-1}_{j,\lambda}\varphi^{2l-1}_{k,\lambda} \right \rVert}^2_2 \frac{j!(2q-1)!}{(j+2q-1)!}\\
    & \hspace{220pt}\frac{k!(2l-1)!}{(k+2l-1)!} \varphi^{2q-1}_{j,\lambda}(2iy',2iv')\varphi^{2l-1}_{k,\lambda}(2iy'',2iv'')\\
   % &=\sum_{l',l''}(2\pi)^{-n}|\mathbf{1}|^n|\la|^{3n}e^{-|\mathbf{1}|\frac{|\la|}{2}(vx-uy)}2^{n/2}(2\pi)^{-2n-2}\sum_{k',k''}{\left \lVert \left (f^{|\mathbf{1}|\la}\times_{|\mathbf{1}||\la|}\varphi^{|\mathbf{1}||\la|,2q}_{l'}\varphi^{|\mathbf{1}||\la|,2l}_{l''}\right)\times_{|\mathbf{1}||\la|}\varphi^{|\mathbf{1}||\la|,2q}_{k'}\varphi^{|\mathbf{1}||\la|,2l}_{k''} \right \rVert}^2_2  \\
    % &  \frac{k'!(2q-1)!}{(k'+2q-1)!} \varphi^{|\mathbf{1}||\la|}_k(2iy',2iv')\frac{k''!(2l-1)!}{(k''+2l-1)!} \varphi^{|\mathbf{1}||\la|,2q}_{k'}
     \\
     &=e^{-\frac{\lambda}{2}(v\cdot x-u\cdot y)}(2\pi)^{-n-3} 3^{-n/2}\sum_{j,k}(\lambda)^{n}C_{j,k}{\left \lVert f^{\boldsymbol{\lambda}}\times_{\lambda}\varphi^{2q-1}_{j,\lambda}\varphi^{2l-1}_{k,\lambda} \right \rVert}^2_2\\
     &\hspace{280pt}\varphi^{2q-1}_{j,\lambda}(2iy',2iv')\varphi^{2l-1}_{k,\lambda}(2iy'',2iv'').
\end{align*}
\end{proof}

\subsection{The case when $m=2$:}
Now we state the Gutzmer formula for reduced H-type group with dimension of center \(m=2\). 
%In this case also, we establish an estimate similar to \eqref{red gutzmer eq3}.
\begin{theorem}\label{Gutzmers formula for reduced}
    Let \(N\) be the H-type group with center \(\mathbb{R}^2\) and \(N_{red}\) be its corresponding reduced H-type group. Let \(f\) be a function on \(N_{red}\) having mean value 0, satisfying the hypothesis of Theorem \ref{theorem 4.2}. Then we have 
    $$
  \begin{aligned}
   &\int_{U \ltimes N_{red}}|f(g.(z, w, \zeta))|^2 d g= \\
&(2 \pi)^{-n } \sum_{\mathbf{p} \neq 0} e^{-2 \mathbf{p} \cdot \tau} e^{-p(v \cdot x-u \cdot y)}\left(\sum_{k=0}^{\infty}\left\|f^{\mathbf{p}} \times_{p} \varphi_{k,p}^{n-1}\right\|_2^2 \frac{k!(n-1)!}{(k+n-1)!} \varphi_{k,p}^{n-1}(-2 i y, -2 iv)\right)|j|^n(2\pi)^{-n-2},
   \end{aligned}
   $$
   where \(|\boldsymbol{p}|=p,\; n=2q\), \(z=x+iy\), \(w=u+iv\) and \(\zeta=\omega+i\tau\).
\end{theorem}
The next proposition is analogue of Proposition \ref{beurling estimate m=3}.
\begin{prop}\label{prop 5.5}
    Let \(N\) be H-type group with center \(\mathbb{R}^2\) and \(f\) be a Schwartz class function on \(N\), Then we have the following estimate
    $$
    \begin{aligned}
        &\displaystyle\int_{U}{\lVert \pi_{\boldsymbol{\lambda}}(\sigma(z,w))\Hat{f}(\boldsymbol{\lambda})\rVert}^2_{HS}d\sigma\\
        &=e^{-\frac{\lambda}{2}(v\cdot x-u\cdot y)}(2\pi)^{-n-2} \sum_{k=0}^{\infty}(\lambda)^{n}{\left \lVert f^{\boldsymbol{\lambda}}\times_{\la}\varphi^{n-1}_{k,\lambda} \right \rVert}^2_2\frac{k!(n-1)!}{(k+n-1)!}\varphi^{n-1}_{k,\lambda}(2iy,2iv),
    \end{aligned}
    $$
    where \(|\boldsymbol{\lambda}|=\lambda,\; n=2q\), \(x=x+iy,w=u+iv\).
\end{prop}
\begin{proof}
The proof follows along the lines of Proposition \ref{beurling estimate m=3}.
\end{proof}
\section{Proof of theorem \ref{Beurling's theorem}}\label{sec6}
We now prove Beurling's theorem.
The Beurling theorem on H-type group \(N \cong \mathbb{C}^n \times \mathbb{R}^m\) with \(m=1,2,3\) is as below. One can see the proof for the case when \(m=1\) in \cite{thangavelubeurling}.
\begin{proof}[Proof of Theorem \ref{Beurling's theorem}]
 Let \(N\) be a H-type group with center \(\mathbb{R}^3\). As in the case of the Heisenberg group, for \(\boldsymbol{\lambda} \in \mathbb{R}^3 \setminus\{0\}\) with \(|\boldsymbol{\lambda}|=\lambda\), consider the function 
$$
\mathcal{F}_{\boldsymbol{\lambda}}(\zeta)=\int_{\mathbb{R}^{2n}}\overline{f^{\boldsymbol{\lambda}}(y,v)}f^{\boldsymbol{\lambda}}(\zeta y, \zeta v) dydv, \;\; \forall \zeta \in \mathbb{C},\; n=2q+2l.
$$
To see that this is a well-defined function, we first show that \(f^{\boldsymbol{\lambda}}\) has a holomorphic extension and the above integral converges.\\
By Proposition \ref{beurling estimate m=3}, we have that 
\begin{align}\label{summation1}
&e^{-\frac{\la}{2}(v\cdot x-u\cdot y)} \sum_{j,k}{\left \lVert f^{\boldsymbol{\lambda}}\times_{\la}\varphi^{2q-1}_{j,\lambda} \varphi^{2l-1}_{k,\lambda}\right \rVert}^2_2\frac{j!(2q-1)!}{(j+2q-1)!}\varphi^{2q-1}_{j,\lambda}(2iy',2iv') \nonumber \\ 
&\hspace{220 pt}\frac{k!(2l-1)!}{(k+2l-1)!}\varphi^{2l-1}_{k,\lambda}(2iy'',2iv'')< \infty.
\end{align}
We recall that \Be \varphi_k^{n-1}={L}^{n-1}_k \left ( \frac{1}{2}(x^2+u^2)\right )e^{-(x^2+u^2)/4}, \text{ where } {L}^{n-1}_k \text { denotes Laguerre polynomial}.\Ee
Hence, for any \(\boldsymbol{\nu} \in \mathbb{R}^3\setminus \{0\}\), 
\begin{align*}
    \varphi^{n-1}_{k,|\boldsymbol{\nu}|}(2iy,2iv)&=\varphi_k^{n-1}(|\boldsymbol{\nu}|^{1/2}(2iy,2iv))\\
    %&=\mathcal{L}^{n-1}_k\left( \frac{1}{2}(-4|\mu|(y^2+v^2)\right)e^{|\mu|(y^2+v^2)}\\
    &={L}^{n-1}_k\left( -2|\boldsymbol{\nu}|(y^2+v^2)\right)e^{|\boldsymbol{\nu}|(y^2+v^2)}\\
    &={L}^{n-1}_k(-2\rho^2)e^{\rho^2}\; \text{ where } \rho^2=|\boldsymbol{\nu}|(y^2+v^2).
\end{align*}
In particular, for \(\boldsymbol{\lambda} \in \mathbb{R}^3 \setminus\{0\}\), \Be \varphi^{2q-1}_{j,\lambda}(2iy',2iv')={L}^{2q-1}_{j}(-2\rho'^2)e^{\rho'^2}, \Ee where \(\rho'^2=|\boldsymbol{\lambda}|(y'^2+v'^2)=\lambda(y'^2+v'^2)\) and \Be \varphi^{2l-1}_{k,\lambda}(2iy'',2iv'')={L}^{2l-1}_{k}(-2\rho''^2)e^{\rho''^2}\Ee where \(\rho''^2=|\boldsymbol{\lambda}|(y''^2+v''^2)=\lambda(y''^2+v''^2)\).
%In particular, for \(\la \in \mathbb{R}^2 \setminus\{0\}\), \(\varphi^{|\mathbf{1}||\la|}_k=\mathcal{L}^{n-1}_k(-2\rho^2)e^{\rho^2}\) where \(\rho^2=|\mathbf{1}||\la|(y^2+v^2)\).
We have the following asymptotic behavior of Laguerre polynomials \cite{szego-orthogonal-polynomials},
$$
{L}^\al_k(s)=\frac{1}{2}\pi^{-1/2}e^{s/2}(-s)^{-\al/2-1/4}e^{2\sqrt{-sk}}(1+\mathcal{O}(k^{-1/2}))
$$
for all \(s\) in the complex plane cut along the positive real axis.\\
Taking \(s=-2\rho'^2\), \(\al=2q-1\), we get 
\Be
{L}^{2q-1}_k(-2\rho'^2)e^{\rho'^2}=\frac{1}{2}\pi^{-\frac{1}{2}}(2\rho'^2)^{-\frac{2q-1}{2}-\frac{1}{4}}k^{\frac{2q-1}{2}-\frac{1}{4}}e^{2\rho'\sqrt{2k}}(1+\mathcal{O}(k^{-1/2})). \Ee
This implies that
\Be \frac{1}{2}\pi^{-\frac{1}{2}}(2\rho'^2)^{-\frac{2q-1}{4}-\frac{1}{4}}e^{2\rho'\sqrt{2k}}k^{\frac{2q-1}{2}-\frac{1}{4}}(-c_1k^{-\frac{1}{2}}+1) \leq {L}^{2q-1}_k(-2\rho'^2)e^{\rho'^2}, \Ee for some $c_1>0$.

By the convergence of the series \eqref{summation1} and the asymptotic behaviour of the Laguerre polynomials, we obtain that for every \(\rho',\rho''>0\),
$$
{\lVert f^{\boldsymbol{\lambda}}\times_{\lambda}\varphi^{2q-1}_{j,\lambda}\varphi^{2l-1}_{k,\lambda}\rVert}^2_2 \leq ce^{-2\rho'\sqrt{(2j+2q)\la}}e^{-2\rho''\sqrt{(2k+2l)\la}}.
$$
Thus, we can write \(f^{\boldsymbol{\lambda}}\) as \(g_{\boldsymbol{\lambda}} \times_{\la}P^{\la}_{\rho',\rho''}\)
 where \(g_{\boldsymbol{\lambda}} \in L^2(\mathbb{R}^{2n})\) and 
$$
P^{\la}_{\rho',\rho''}=(2\pi)^{-2q-2l}\sum_{j,k}e^{-\rho'\sqrt{(2j+2q)\la}}\varphi^{2q-1}_{j,\lambda}e^{-\rho''\sqrt{(2k+2l)\la}}\varphi^{2l-1}_{k,\lambda}.
$$
For the case when \(m=2\), by Proposition \ref{prop 5.5}, we have that 
\begin{equation}\label{summation2}
e^{-\frac{\la}{2}(vx-uy)} \sum_{k=0}^{\infty}{\left \lVert f^{\boldsymbol{\lambda}}\times_{\la}\varphi^{2q-1}_{k,\lambda} \right \rVert}^2_2\frac{k!(2q-1)!}{(k+2q-1)!}\varphi^{2q-1}_{k,\lambda}(2iy,2iv)< \infty.
\end{equation}
Proceeding as above, the convergence of the series \eqref{summation2} implies that for every \(\rho >0\), 
$$
{\lVert f^{\boldsymbol{\lambda}}\times_{\la}\varphi^{2q-1}_{k,\lambda}\rVert}^2_2 \leq ce^{-2\rho\sqrt{(2k+2q)\la}}.
$$
Hence we can write the function \(f^{\boldsymbol{\lambda}}\) as \(g_{\boldsymbol{\la}}\times_{\la}P^{\la}_\rho\) where \(g_{\la} \in L^2(\mathbb{R}^{2n})\), \(n=2q\) and 
$$
P_{\rho}^{\la}=(2\pi)^{-2q}\sum_{k=0}^{\infty}e^{-\rho\sqrt{(2k+2q)\la}}\varphi^{2q-1}_{k,\lambda}.
$$
%The rest of the proof is verbatim to proof of \cite[Theorem 1.2]{thangavelubeurling}. 
This way we can extend the function \(f^{\boldsymbol{\lambda}}\) holomorphically to whole of \(\mathbb{C}^{2n}\). Following the similar steps as in \cite[Theorem 1.2]{thangavelubeurling}, we find that the function \( \mathcal{F}_{\boldsymbol{\lambda}} \) defines a holomorphic function. Subsequently, following the proof of Hedenmalm \cite{Hedenmalm} we can derive the desired result.
\end{proof}

%\section{Radon Transform on H-type groups}\label{sec6}
%In this section we recall a few terminologies related to Radon transform on H-type groups and prove two more versions of the Beurling's theorem. As before, throughout this section \(N\) refers to H-type group unless stated otherwise. To use Radon transform, we need to use certain different unitary representations of H-type groups that come from the unitary representations of the Heisenberg group. 

%\begin{remark}
 %   Note that while establishing Theorem \ref{5.4}, we had to restrict ourselves to the case where the dimension of the center \(m=1,2,3\), whereas, in Theorems \ref{first Beurling's using radon transform} and \ref{second Beurling's using radon transform}  there is no restriction on the dimension of the center.
%\end{remark}

  \bibliographystyle{acm}
\bibliography{ref}

\begin{thebibliography}{10}

\bibitem{bonfi}
{\sc Bonfiglioli, A., Lanconelli, E., and Uguzzoni, F.}
\newblock {\em Stratified {L}ie groups and potential theory for their sub-{L}aplacians}.
\newblock Springer Monographs in Mathematics. Springer, Berlin, 2007.

\bibitem{dasgupta-santosh-localization-operator}
{\sc Dasgupta, A., and Nayak, S.~K.}
\newblock Localization operator and {W}eyl transform on reduced {H}eisenberg group with multidimensional center.
\newblock {\em J. Math. Phys. 64}, 12 (2023), Paper No. 123506, 20.

\bibitem{dasgupta-thangavelu-nilmanifold}
{\sc Dasgupta, A., and Thangavelu, S.}
\newblock Heat kernel transform on nilmanifolds associated to {H}-type groups.
\newblock {\em Tohoku Math. J. (2) 64}, 3 (2012), 439--451.

\bibitem{dasguptawong}
{\sc Dasgupta, A., and Wong, M.~W.}
\newblock Weyl transforms for h-type groups.
\newblock {\em Journal of Pseudo-Differential Operators and Applications 6}, 1 (2015), 11--19.

\bibitem{Hedenmalm}
{\sc Hedenmalm, H.}
\newblock Heisenberg's uncertainty principle in the sense of {B}eurling.
\newblock {\em J. Anal. Math. 118}, 2 (2012), 691--702.

\bibitem{kaplanricci}
{\sc Kaplan, A., and Ricci, F.}
\newblock Harmonic analysis on groups of {H}eisenberg type.
\newblock In {\em Harmonic analysis ({C}ortona, 1982)}, vol.~992 of {\em Lecture Notes in Math.} Springer, Berlin, 1983, pp.~416--435.

\bibitem{functionalhtype}
{\sc Liu, H., and Song, M.}
\newblock A functional calculus and restriction theorem on h-type groups.
\newblock {\em Pacific Journal of Mathematics 286}, 2 (2017), 291--305.

\bibitem{LiurestrictionHtype}
{\sc Liu, H., and Wang, Y.}
\newblock A restriction theorem for the {H}-type groups.
\newblock {\em Proc. Amer. Math. Soc. 139}, 8 (2011), 2713--2720.

\bibitem{narayananinjectivity}
{\sc Narayanan, E.~K., Sanjay, P.~K., and Yasser, K.~T.}
\newblock Injectivity of spherical means on {$H$}-type groups.
\newblock {\em J. Fourier Anal. Appl. 30}, 3 (2024), Paper No. 33, 24.

\bibitem{Beurling-rudra-parui}
{\sc Parui, S., and Sarkar, R.~P.}
\newblock Beurling's theorem and {$L^p$}-{$L^q$} {M}organ's theorem for step two nilpotent {L}ie groups.
\newblock {\em Publ. Res. Inst. Math. Sci. 44}, 4 (2008), 1027--1056.

\bibitem{riccicomm}
{\sc Ricci, F.}
\newblock Commutative algebras of invariant functions on groups of {H}eisenberg type.
\newblock {\em J. London Math. Soc. (2) 32}, 2 (1985), 265--271.

\bibitem{roncal-thangavelu-holomorphic-extension}
{\sc Roncal, L., and Thangavelu, S.}
\newblock Holomorphic extensions of eigenfunctions on {$NA$} groups.
\newblock {\em arXiv preprint arXiv:2005.09894\/} (2020).

\bibitem{sarkarbeurling}
{\sc Sarkar, R.~P., and Sengupta, J.}
\newblock Beurling's theorem for {R}iemannian symmetric spaces. {II}.
\newblock {\em Proc. Amer. Math. Soc. 136}, 5 (2008), 1841--1853.

\bibitem{strichartz}
{\sc Strichartz, R.~S.}
\newblock {$L^p$} harmonic analysis and {R}adon transforms on the {H}eisenberg group.
\newblock {\em J. Funct. Anal. 96}, 2 (1991), 350--406.

\bibitem{szego-orthogonal-polynomials}
{\sc Szeg\H~o, G.}
\newblock {\em Orthogonal polynomials}, third~ed., vol.~Vol. 23 of {\em American Mathematical Society Colloquium Publications}.
\newblock American Mathematical Society, Providence, RI, 1967.

\bibitem{gutzmer}
{\sc Thangavelu, S.}
\newblock Gutzmer's formula and {P}oisson integrals on the {H}eisenberg group.
\newblock {\em Pacific J. Math. 231}, 1 (2007), 217--237.

\bibitem{thangaveluheisenberg}
{\sc Thangavelu, S.}
\newblock {\em Harmonic analysis on the Heisenberg group}, vol.~159.
\newblock Springer Science \& Business Media, 2012.

\bibitem{thangavelubeurling}
{\sc Thangavelu, S.}
\newblock Beurling's theorem on the {H}eisenberg group.
\newblock {\em Ark. Mat. 60}, 2 (2022), 417--442.

\bibitem{yangheatkernel}
{\sc Yang, Q., and Zhu, F.}
\newblock The heat kernel on {H}-type groups.
\newblock {\em Proc. Amer. Math. Soc. 136}, 4 (2008), 1457--1464.

\bibitem{yasser2024gelfand}
{\sc Yasser, K.}
\newblock Gelfand pairs and spherical means on h-type groups.
\newblock {\em Indian Journal of Pure and Applied Mathematics 55}, 4 (2024), 1316--1328.

\bibitem{mingkai}
{\sc Yin, M., and He, J.}
\newblock H-type group, a-weyl transform and pseudo-differential operators.
\newblock {\em International Journal of Mathematical Analysis 9}, 25 (2015), 1201--1214.

\end{thebibliography}

\end{document}